\UseRawInputEncoding
\documentclass[12pt,twoside, epsf]{article}

\usepackage{color,graphicx,times}

\makeatother

\let \ttorg \tt \def \tt{\ttorg \obeyspaces}

\begin{document}

\newcommand{\Across}{\raisebox{-0.25\height}{\includegraphics[width=0.5cm]{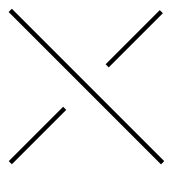}}}
\newcommand{\Asmooth}{\raisebox{-0.25\height}{\includegraphics[width=0.5cm]{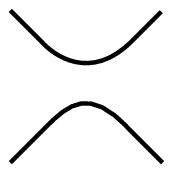}}}
\newcommand{\Bsmooth}{\raisebox{-0.25\height}{\includegraphics[width=0.5cm]{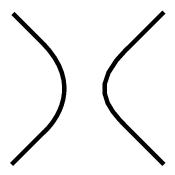}}}
\newcommand{\Rcurl}{\raisebox{-0.25\height}{\includegraphics[width=0.5cm]{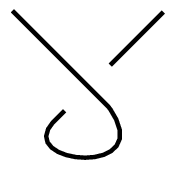}}}
\newcommand{\Lcurl}{\raisebox{-0.25\height}{\includegraphics[width=0.5cm]{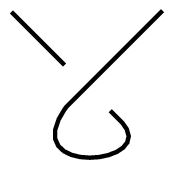}}}
\newcommand{\Arc}{\raisebox{-0.25\height}{\includegraphics[width=0.5cm]{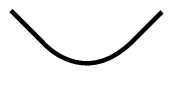}}}

\title{\Large\bf An Introduction to Khovanov Homology}
\author{Louis H. Kauffman\\ Department of Mathematics, Statistics \\ and Computer Science (m/c
249)    \\ 851 South Morgan Street   \\ University of Illinois at Chicago\\
Chicago, Illinois 60607-7045\\ $<$kauffman@uic.edu$>$\\}

\maketitle

\thispagestyle{empty}

\subsection*{\centering Abstract}

{\em This paper is an introduction to Khovanov Homology.}
\bigbreak

\noindent {\bf Keywords.} bracket polynomial, Khovanov homology, cube category, simplicial category, tangle cobordism, chain complex, chain homotopy, unitary transformation, quantum computing, quantum information theory, link homology, categorification.

\section{Introduction}
This paper is an introduction to Khovanov homology.\\

We start with a quick introduction to the bracket polynomial, reformulating it and the Jones polynomial
so that the value of an unknotted loop is $q + q^{-1}.$ We then introduce enhanced states for the 
bracket state sum so that, in terms of these enhanced states the bracket is a sum of monomials.\\

In Section 3, we point out that the shape of the collection of bracket states for a given diagram is a cube and that this cube can be taken to be a category. It is an example of a {\em cube category}. We show that functors from a cube category to a category of modules naturally have homology theories associated with them. In Section 4 we show how to make a homology theory (Khovanov homology) from the states of the bracket so that the enhanced states are the generators of the chain complex. We show how a Frobenius algebra structure arises naturally from this adjacency structure for the enhanced states.
Finally we show that the resulting homology is an example of homology related to a module functor on the cube category as described in Section 3.\\

 In Section 5, we give a short exposition of Dror BarNatan's tangle cobordism theory for Khovanov homology. This theory replaces Khovanov homology by an abstract chain homotopy class of a complex of surface cobordisms associated with the states of a knot or link diagram. Topologically equivalent links give rise to abstract chain complexes (special cobordism categories) that are chain homotopy equivalent. We show how the $4Tu$ tubing relation of BarNatan is exactly what is needed to show that these chain complexes are invariant under the second Reidemeister move. This is the key ingredient in the full invariance under Reidemeister moves, and it shows how one can reinvent the $4Tu$ relation by searching for that homotopy. Once one has the $4Tu$ relation it is easy to see that it
 is equivalent to the tube-cutting relation that is satisfied by the Frobenius algebra we have already discussed. In this way we obtain both the structure of the abstract chain homotopy category, and the invariance of Khovanov homology under the Reidemeister moves.
 \\
 
 In Section 6 we show how the tube-cutting relation can be used to derive a class of Frobenius algebras depending on a choice of parameter $t$ in the base field.
When $t=0$ we have the original Frobenius algebra for Khovanov homology. For $t=1$ we have the Lee Algebra on which is based the Rasmussen invariant.
The derivation in this section of a class of Frobenius algebras from the tube-cutting relation, shows that one can begin Khovanov homology with the abstract categorical chain complex associated with the 
Cube Category of a link and from this data find the Frobenius algebras that can produce the actual homology theories.\\

In Section 7 we give a short exposition of the Rasmussen invariant and its application to 
finding the four-ball genus of torus knots.\\

 In Section 8 we give a description of Khovanov homology as the homology of a simplicial module by following our description of the cube category in this context.
In Section 9 we discuss a quantum context for Khovanov homology that is obtained by building a Hilbert space whose orthonormal basis is the 
set of enhanced states of a diagram $K.$ Then there is a unitary transformation $U_{K}$  of this Hilbert space so that the Jones polynomial $J_{K}$ is the trace of $U_{K}: J_{K} = Trace(U_{K}).$ We discuss a generalization where the linear space of the Khovanov homology itself is taken to be the Hilbert space.
In this case we can define a unitary transformation $U'_{K}$ so that, for values of $q$ and $t$ on the unit 
circle, the Poincar\'e polynomial for the Khovanov homology is the trace of $U'_{K}.$ Section 10 is
is a discussion,with selected references, of other forms of link homology and categorification, including generalizations of Khovanov homology to virtual knot theory.
\bigbreak

It gives the author great pleasure to thank the members of the Quantum Topology Seminar at the University of Illinois at Chicago for many useful conversations and to thank the Perimeter Institute
in Waterloo, Canada for their hospitality while an early version of this  paper was being completed. The present paper is an extension of the paper \cite{KTrieste}.
\bigbreak

\section{Bracket Polynomial and Jones Polynomial}
The bracket polynomial \cite{KaB} model for the Jones polynomial \cite{JO,JO1,JO2,Witten} is usually described by the expansion
$$\langle \Across \rangle=A \langle \Asmooth \rangle + A^{-1}\langle \Bsmooth \rangle$$
Here the small diagrams indicate parts of otherwise identical larger knot or link diagrams. The two types of smoothing (local diagram with no crossing) in
this formula are said to be of type $A$ ($A$ above) and type $B$ ($A^{-1}$ above).

$$\langle \bigcirc \rangle = -A^{2} -A^{-2}$$
$$\langle K \, \bigcirc \rangle=(-A^{2} -A^{-2}) \langle K \rangle $$
$$\langle \Rcurl \rangle=(-A^{3}) \langle \Arc \rangle $$
$$\langle \Lcurl \rangle=(-A^{-3}) \langle \Arc \rangle $$
One uses these equations to normalize the invariant and make a model of the Jones polynomial.
In the normalized version we define $$f_{K}(A) = (-A^{3})^{-wr(K)} \langle K \rangle / \langle \bigcirc \rangle $$
where the writhe $wr(K)$ is the sum of the oriented crossing signs for a choice of orientation of the link $K.$ Since we shall not use oriented links
in this paper, we refer the reader to \cite{KaB} for the details about the writhe. One then has that $f_{K}(A)$ is invariant under the Reidemeister moves
(again see \cite{KaB}) and the original Jones polynomial $V_{K}(t)$ is given by the formula $$V_{K}(t) = f_{K}(t^{-1/4}).$$ The Jones polynomial
has been of great interest since its discovery in 1983 due to its relationships with statistical mechanics, due to its ability to often detect the
difference between a knot and its mirror image and due to the many open problems and relationships of this invariant with other aspects of low
dimensional topology.
\bigbreak

\noindent {\bf The State Summation.} In order to obtain a closed formula for the bracket, we now describe it as a state summation.
Let $K$ be any unoriented link diagram. Define a {\em state}, $S$, of $K$  to be the collection of planar loops resulting from  a choice of
smoothing for each  crossing of $K.$ There are two choices ($A$ and $B$) for smoothing a given  crossing, and
thus there are $2^{c(K)}$ states of a diagram with $c(K)$ crossings.
In a state we label each smoothing with $A$ or $A^{-1}$ according to the convention
indicated by the expansion formula for the bracket. These labels are the  {\em vertex weights} of the state.
There are two evaluations related to a state. The first is the product of the vertex weights,
denoted $\langle K|S \rangle .$
The second is the number of loops in the state $S$, denoted  $||S||.$

\noindent Define the {\em state summation}, $\langle K \rangle $, by the formula

$$\langle K \rangle  \, = \sum_{S} <K|S> \delta^{||S||}$$
where $\delta = -A^{2} - A^{-2}.$
This is the state expansion of the bracket. It is possible to rewrite this expansion in other ways. For our purposes in
this paper it is more convenient to think of the loop evaluation as a sum of {\it two} loop evaluations, one giving $-A^{2}$ and one giving
$-A^{-2}.$ This can be accomplished by letting each state curve carry an extra label of $+1$ or $-1.$ We describe these {\it enhanced states} below. But before we do this, it will be useful for the reader to examine Figure 2. In Figure 2 we show all the states for the right-handed trefoil knot, labelling the sites
with $A$ or $B$ where $B$ denotes a smoothing that would receive $A^{-1}$ in the state expansion.
\bigbreak

Note that in the state enumeration in Figure 2 we have organized the states in tiers so that the state
that has only $A$-smoothings is at the top and the state that has only $B$-smoothings is at the bottom.
\bigbreak

\begin{figure}
     \begin{center}
     \begin{tabular}{c}
     \includegraphics[width=6cm]{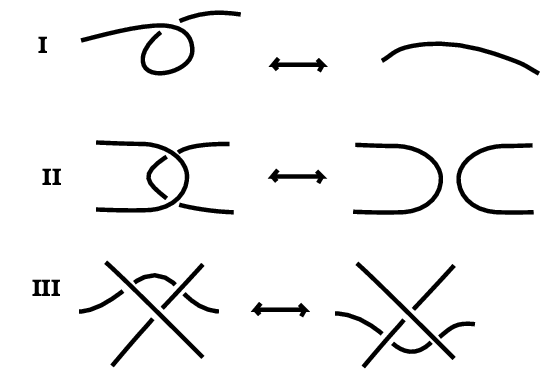}
     \end{tabular}
     \caption{\bf Reidemeister Moves}
     \label{Figure 1}
\end{center}
\end{figure}

\noindent {\bf Changing Variables.} Letting $c(K)$ denote the number of crossings in the diagram $K,$ if we replace $\langle K
\rangle$ by
$A^{-c(K)} \langle K \rangle,$ and then replace $A^2$ by $-q^{-1},$ the bracket is then rewritten in the
following form:
$$\langle \Across \rangle=\langle \Asmooth \rangle-q\langle \Bsmooth \rangle $$
with $\langle \bigcirc\rangle=(q+q^{-1})$.
It is useful to use this form of the bracket state sum
for the sake of the grading in the Khovanov homology (to be described below). We shall
continue to refer to the smoothings labeled $q$ (or $A^{-1}$ in the
original bracket formulation) as {\it $B$-smoothings}.
\bigbreak

We catalog here the resulting behaviour of this modified bracket under the Reidemeister moves.
$$\langle \bigcirc \rangle = q + q^{-1}$$
$$\langle K \, \bigcirc \rangle=(q + q^{-1}) \langle K \rangle $$
$$\langle \Rcurl \rangle=q^{-1} \langle \Arc \rangle $$
$$\langle \Lcurl \rangle= -q^{2} \langle \Arc \rangle $$
$$\langle \raisebox{-0.25\height}{\includegraphics[width=0.5cm]{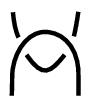}}\rangle =
-q \langle  \raisebox{-0.50\height}{\includegraphics[width=0.5cm]{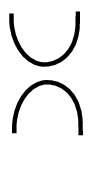}}\rangle$$
$$\langle \raisebox{-0.25\height}{\includegraphics[width=0.8cm]{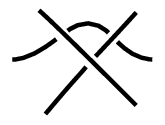}}\rangle =
\langle  \raisebox{-0.25\height}{\includegraphics[width=0.8cm]{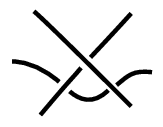}}\rangle$$

It follows that if we define $$J_{K} = (-1)^{n_{-}} q^{n_{+} - 2n_{-}} \langle K \rangle,$$
where $n_{-}$ denotes the number of negative crossings in $K$ and $n_{+}$ denotes the number
of positive crossings in $K$, then $J_{K}$ is invariant under all three Reidemeister moves.
Thus $J_{K}$ is a version of the Jones polynomial taking the value $q + q^{-1}$ on an unknotted circle.
\bigbreak

\noindent {\bf Using Enhanced States.}
We now use the convention of {\it enhanced
states} where an enhanced state has a label of $1$ or $-1$ on each of
its component loops. We then regard the value of the loop $q + q^{-1}$ as
the sum of the value of a circle labeled with a $1$ (the value is
$q$) added to the value of a circle labeled with an $-1$ (the value
is $q^{-1}).$ We could have chosen the less neutral labels of $+1$ and $X$ so that
$$q^{+1} \Longleftrightarrow +1 \Longleftrightarrow 1$$
and
$$q^{-1} \Longleftrightarrow -1 \Longleftrightarrow x,$$
since an algebra involving $1$ and $x$ naturally appears later in relation to Khovanov homology. It does no harm to take this form of labeling from the
beginning. The use of enhanced states for formulating Khovanov homology was pointed out by Oleg Viro in
\cite{Viro}.
\bigbreak

Consider the form of the expansion of this version of the
bracket polynonmial in enhanced states. We have the formula as a sum over enhanced states $s:$
$$\langle K \rangle = \sum_{s} (-1)^{i(s)} q^{j(s)} $$
where $i(s)$ is the number of $B$-type smoothings in $s$ and $j(s) = i(s) + \lambda(s)$, with $\lambda(s)$ the number of loops  labeled $1$ minus
the number of loops labeled $-1$ in the enhanced state $s.$
\bigbreak

One advantage of the expression of the bracket polynomial via enhanced states is that it is now a sum of monomials. We shall make use of this property
throughout the rest of the paper.
\bigbreak

\section{Khovanov Homology and the Cube Category}
We are going to make a chain complex from the states of the bracket polynomial so that the homology
of this chain complex is a knot invariant. One way to see how such a homology theory arises is to step back and note that the collection of states for a diagram $K$ forms a category in the shape of a cube.
A functor from such a category to a category of modules gives rise to a homology theory in a natural way, as we explain below.
\bigbreak

Examine Figure 2 and Figure 3. In Figure 2 we show all the 
standard bracket states for the trefoil knot with arrows between them whenever the state at the output of the arrow is obtained from the state at the input of the arrow by a single smoothing of a site of type $A$ to a site of type $B$. The abstract structure of this collection of states is a category with objects of the form
$\langle ABA \rangle $ where this symbol denotes one of the states in the state diagram of Figure 2. In Figure 3 we illustrate this  cube  category (the states are arranged in the form of a cube)
by replacing the states in Figure 2 by symbols $\langle XYZ \rangle$ where each literal is either an $A$ or a $B.$
A typical generating morphism in the $3$-cube category is $$\langle ABA\rangle \longrightarrow \langle BBA\rangle .$$
\bigbreak

\begin{figure}
     \begin{center}
     \begin{tabular}{c}
     \includegraphics[width=7cm]{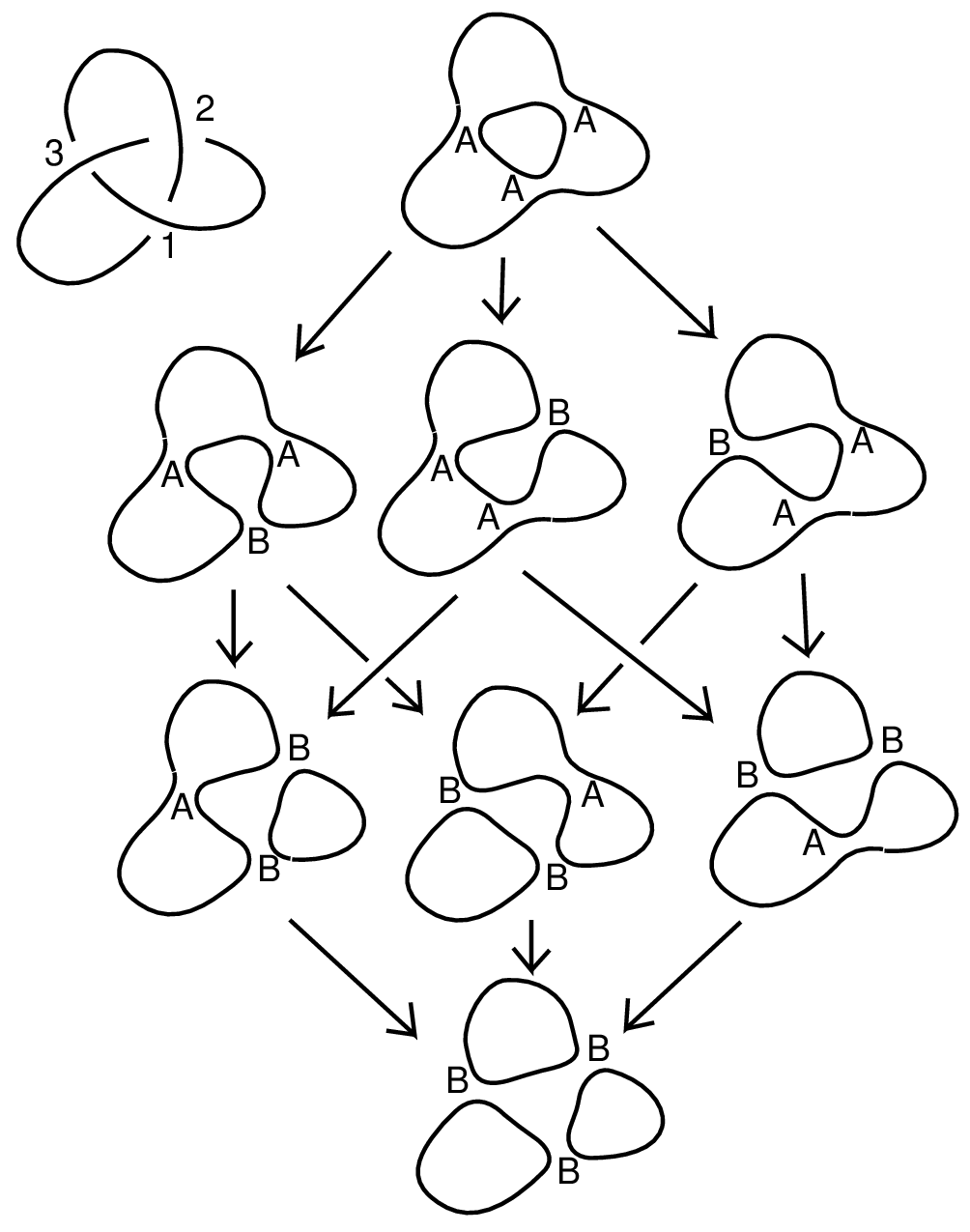}
     \end{tabular}
     \caption{\bf Bracket States and Khovanov Complex}
    \label{Figure 2}
\end{center}
\end{figure}

\begin{figure}
     \begin{center}
     \begin{tabular}{c}
     \includegraphics[width=7cm]{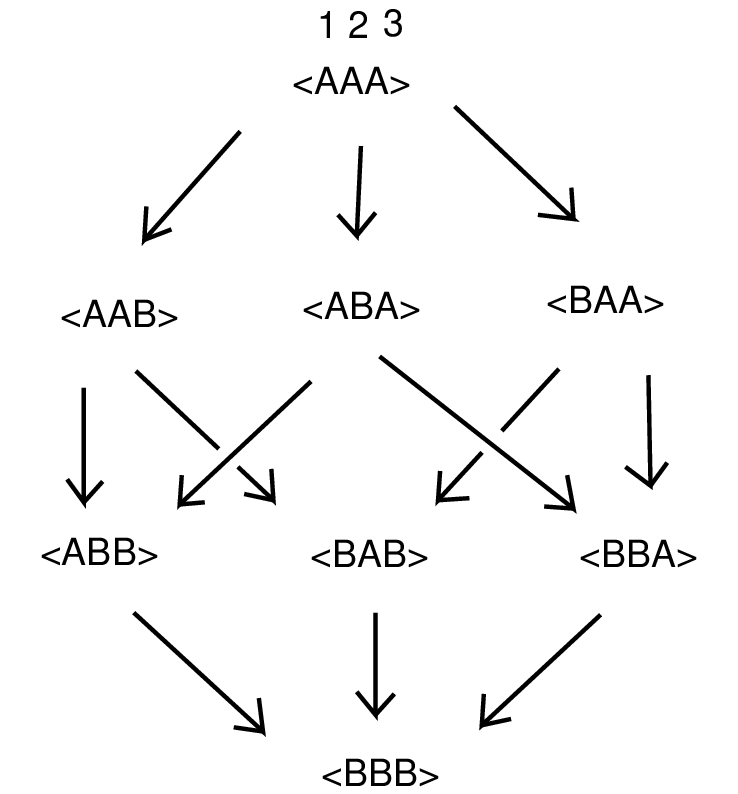}
     \end{tabular}
     \caption{\bf Cube Category}
    \label{Figure 3}
\end{center}
\end{figure}

We formalize this way of looking at the bracket states as follows. Let $\mathcal{S}(K)$ denote a category 
associated with the states of the bracket for a diagram $K$ whose objects are the states, with sites
labeled $A$ and $B$ as in Figure 2. A morphism in this category is an arrow from a state with a given
number of $A$'s to a state with fewer $A$'s.
\bigbreak

\noindent Let ${\mathcal D}^{n} = \{A,B \}^{n}$ be the $n$-cube category whose objects are the $n$-sequences from
the set $\{A, B \}$ and whose morphisms are arrows from sequences with greater numbers of $A$'s
to sequences with fewer numbers of $A$'s. Thus ${\mathcal D}^{n}$ is equivalent to the poset category
of subsets of $\{ 1,2,\cdots n \}.$ We make a functor $\mathcal{R}: {\mathcal D}^{n} \longrightarrow \mathcal{S}(K)$
for a diagram $K$ with $n$ crossings as follows.
We map sequences in the cube category  to bracket states
by choosing to label the crossings of the diagram $K$ from the set  $\{ 1,2,\cdots n \},$ and letting
this  functor take abstract $A$'s and $B$'s in the cube category to smoothings at those crossings of type $A$ or type $B.$ Thus each sequence in the cube category is associated with a unique state of $K$
when $K$ has $n$ crossings. By the same token, we define a functor 
$\mathcal{S}:\mathcal{S}(K) \longrightarrow {\mathcal D}^{n} $ by associating a sequence to each state
and morphisms between sequences corresponding to the state smoothings. With these conventions,
the two compositions of these morphisms are the identity maps on their respective categories.
\bigbreak

Let ${\mathcal M }$ be a pointed category with finite sums, and let
${\mathcal F} : {\mathcal D}^{n}  \longrightarrow  {\mathcal M } $ be a functor.  
In our case ${\mathcal M }$ will be a category of
modules and ${\mathcal F}$ will carry $n$-sequences to certain tensor powers corresponding to the standard bracket states of a knot or link $K.$ We postpone this construction for a moment, and point out
that there is a natural structure of chain complex associated with the functor ${\mathcal F}.$
First note that each object in ${\mathcal D}^{n}$ has the form 
$$X = \langle X_{0} \cdots X_{n-1} \rangle$$
where each $X_{i}$ equals either $A$ or $B$ and we have morphisms 
$$d_{i}: \langle X_{0} \cdots X_{i} \cdots X_{n-1} \rangle \longrightarrow 
\langle X_{0} \cdots \bar{X}_{i} \cdots X_{n-1} \rangle$$ 
whenever $X_{i} = A$ and (by definition)  $\bar{X}_{i} = B.$  We then define
$$\partial_{i} = {\mathcal C}(d_{i}): {\mathcal C}\langle X_{0} \cdots X_{i} \cdots X_{n-1} \rangle
\longrightarrow {\mathcal C}\langle X_{0} \cdots \bar{X}_{i} \cdots X_{n-1} \rangle$$
whenever $d_{i}$ is defined. We then define the chain complex$C$ by 
$$C^{k} = \bigoplus_{X} {\mathcal C}\langle X_{0} \cdots  X_{n-1} \rangle $$
where each sequence $X = \langle X_{0} \cdots  X_{n-1} \rangle$ has $k$
$B$'s. With this we define $$\partial: C^{k} \longrightarrow C^{k+1}$$ by the formula
$$\partial x = \Sigma_{i=0}^{n-1} (-1)^{c(X,i)} \partial_{i}(x)$$ for
 $x \in {\mathcal C}X = {\mathcal C}\langle X_{0} \cdots  X_{n-1} \rangle $
and $c(X,i)$ denotes the number of $A$'s in the sequence $X$ that precede $X_{i}.$
\bigbreak

We want $\partial^{2} = 0$ and it is easy to see that this is equivalent to the condition that
$\partial_{i} \partial_{j} = \partial_{j}  \partial_{i}$ for $i \neq j$ whenever these maps and compositions are defined. We can assume that the functor ${\mathcal F}$ has this property, or we can build it in axiomatically by adding the corresponding relations to the cube category in the form 
$d_{i}  d_{j} = d_{j}  d_{i}$ for $i \neq j$ whenever these maps are defined. In the next section we shall see that there is a natural way to define the maps in the state category so that this
condition holds. Once we axiomatize this commutation relation at the level of the state category or the cube category, then the functor ${\mathcal F}$ will induce a chain complex and homology as above.
\bigbreak

In this way, we see that a suitable functor from the cube category to a module category allows  us to define homology
that is modeled on the ``shape" of the cube. The set of bracket states forms
a natural functorial image of the cube category, and that makes it possible to define the Khovanov chain complex. In terms of the bracket states, we will map each state loop to a specific module $V$, and each
state to a tensor power of $V$ to the number of loops in the state. The details of this construction are in the next section. 
\bigbreak

We use a specific construction for the Khovanov complex that is directly related to the enhanced states for the bracket polynomial, as we will see in the next section. In this construction
we will use the enhanced states, regarding each loop as labeled with either $1$ or $x$ for a module
$V = k[x]/(x^{2})$ associated with the loop (where $k = Z/2Z$ or $k = Z$.) Thus the two labelings of the loop will corespond to the two generators of the module $V.$ A state that is a collection of loops will be associated with $V^{\otimes r}$ where $r$ is the number of loops in the state. In this way we will obtain a functor from the state category
to a module category, and at the same time it will happen that any single enhanced state will correspond to a generator of the chain complex. In the next section we show how naturally this algebra appears in 
relation to the enhanced states. We then return to the categorical point of view and see how, surface cobordisms of circles provide an abstract category for the invariant. 
\bigbreak

\section{Khovanov Homology}

In this section, we describe Khovanov homology
along the lines of \cite{Kho,BN}, and we tell the story so that the gradings and the structure of the differential emerge in a natural way.
This approach to motivating the Khovanov homology uses elements of Khovanov's original approach, Viro's use of enhanced states for the bracket
polynomial \cite{Viro}, and Bar-Natan's emphasis on tangle cobordisms \cite{BN0,BN}. We use similar considerations in our paper \cite{DKM}.
\bigbreak

Two key motivating ideas are involved in finding the Khovanov invariant. First of all, one would like to {\it categorify} a link polynomial such as
$\langle K \rangle.$ There are many meanings to the term categorify, but here the quest is to find a way to express the link polynomial
as a {\it graded Euler characteristic} $\langle K \rangle = \chi_{q} \langle {\mathcal H}(K) \rangle$ for some homology theory associated with $\langle K \rangle.$
\bigbreak

We will use the bracket polynomial and its enhanced states as described in the previous sections of this paper.
To see how the Khovanov grading arises, consider the form of the expansion of this version of the
bracket polynomial in enhanced states. We have the formula as a sum over enhanced states $s:$
$$\langle K \rangle = \sum_{s} (-1)^{i(s)} q^{j(s)} $$
where $i(s)$ is the number of $B$-type smoothings in $s$, $\lambda(s)$ is the number of loops in $s$ labeled $1$ minus the number of loops
labeled $X,$ and $j(s) = i(s) + \lambda(s)$.
This can be rewritten in the following form:
$$\langle K \rangle  =  \sum_{i \,,j} (-1)^{i} q^{j} dim({\mathcal C}^{ij}) $$
where we define ${\mathcal C}^{ij}$ to be the linear span (over the complex numbers for the purpose of this paper, but over the integers or the integers modulo
two for other contexts) of the set of enhanced states with
$i(s) = i$ and $j(s) = j.$ Then the number of such states is the dimension $dim({\mathcal C}^{ij}).$
\bigbreak

\noindent We would like to have a  bigraded complex composed of the ${\mathcal C}^{ij}$ with a
differential
$$\partial:{\mathcal C}^{ij} \longrightarrow {\mathcal C}^{i+1 \, j}.$$
The differential should increase the {\it homological grading} $i$ by $1$ and preserve the
{\it quantum grading} $j.$
Then we could write
$$\langle K \rangle = \sum_{j} q^{j} \sum_{i} (-1)^{i} dim({\mathcal C}^{ij}) = \sum_{j} q^{j} \chi({\mathcal C}^{\bullet \, j}),$$
where $\chi({\mathcal C}^{\bullet \, j})$ is the Euler characteristic of the subcomplex ${\mathcal C}^{\bullet \, j}$ for a fixed value of $j.$
\bigbreak

\noindent This formula would constitute a categorification of the bracket polynomial. Below, we
shall see how {\it the original Khovanov differential $\partial$ is uniquely determined by the restriction that $j(\partial s) = j(s)$ for each enhanced state
$s$.} Since $j$ is
preserved by the differential, these subcomplexes ${\mathcal C}^{\bullet \, j}$ have their own Euler characteristics and homology. We have
$$\chi(H({\mathcal C}^{\bullet \, j})) = \chi({\mathcal C}^{\bullet \, j}) $$ where $H({\mathcal C}^{\bullet \, j})$ denotes the homology of the complex
${\mathcal C}^{\bullet \, j}$. We can write
$$\langle K \rangle = \sum_{j} q^{j} \chi(H({\mathcal C}^{\bullet \, j})).$$
The last formula expresses the bracket polynomial as a {\it graded Euler characteristic} of a homology theory associated with the enhanced states
of the bracket state summation. This is the categorification of the bracket polynomial. Khovanov proves that this homology theory is an invariant
of knots and links (via the Reidemeister moves of Figure 1), creating a new and stronger invariant than the original Jones polynomial.
\bigbreak

We will construct the differential in this complex first for mod-$2$ coefficients.
The differential is based on regarding two states as {\it adjacent} if one differs from the other by a single smoothing at some site.
Thus if $(s,\tau)$ denotes a pair consisting in an enhanced state $s$ and site $\tau$ of that state with $\tau$ of type $A$, then we consider
all enhanced states $s'$ obtained from $s$ by smoothing at $\tau$ and relabeling only those loops that are affected by the resmoothing.
Call this set of enhanced states $S'[s,\tau].$ Then we shall define the {\it partial differential} $\partial_{\tau}(s)$ as a sum over certain elements in
$S'[s,\tau],$ and the differential by the formula $$\partial(s) = \sum_{\tau} \partial_{\tau}(s)$$ with the sum over all type $A$ sites $\tau$ in $s.$
It then remains to see what are the possibilities for $\partial_{\tau}(s)$ so that $j(s)$ is preserved.
\bigbreak

\noindent Note that if $s' \in S'[s,\tau]$, then $i(s') = i(s) + 1.$ Thus $$j(s') = i(s') + \lambda(s') = 1 + i(s) + \lambda(s').$$ From this
we conclude that $j(s) = j(s')$ if and only if $\lambda(s') = \lambda(s) - 1.$ Recall that
$$\lambda(s) = [s:+] - [s:-]$$ where $[s:+]$ is the number of loops in $s$ labeled $+1,$ $[s:-]$ is the number of loops
labeled $-1$ (same as labeled with $x$) and $j(s) = i(s) + \lambda(s)$.
\bigbreak

In the following proposition we assume that the partial derivatives $\partial_{\tau}(s)$ are local in the sense that the loops that are not affected by the resmoothing are not relabeled (just as we have indicated in the previous paragraph). We also assume that the maps we define for partial differentials do not vanish unless this is forced by the grading, and that coefficients of individual tensor products
are taken to be equal to $1.$ In other words, we see that if we take the ``simplest" partial differentials that leave  $j(s)$ invariant, then the differentials are determined by this condition. It is interesting to see how this
works. We shall see later in the paper that there are deeper and more elegant ways to find the algebra indicated below.\\

\noindent {\bf Proposition.} The partial differentials $\partial_{\tau}(s)$ are determined (in the above sense) by the condition that $j(s') = j(s)$ for all $s'$
involved in the action of the partial differential on the enhanced state $s.$ This form of the partial differential can be described by the
following structures of multiplication and comultiplication on the algebra $V = k[x]/(x^{2})$ where
 $k = Z/2Z$ for mod-2 coefficients, or $k = Z$
for integral coefficients.
\begin{enumerate}
\item The element $1$ is a multiplicative unit and $x^2 = 0.$
\item $\Delta(1) = 1 \otimes x + x \otimes 1$ and $\Delta(x) = x \otimes x.$
\end{enumerate}
These rules describe the local relabeling process for loops in a state. Multiplication corresponds to the case where two loops merge to a single loop,
while comultiplication corresponds to the case where one loop bifurcates into two loops.
\bigbreak

\noindent {\bf Proof.}
Using the above description of the differential, suppose that
there are two loops at $\tau$ that merge in the smoothing. If both loops are labeled $1$ in $s$ then the local contribution to $\lambda(s)$ is $2.$
Let $s'$ denote a smoothing in $S[s,\tau].$ In order for the local $\lambda$ contribution to become $1$, we see that the merged loop must be labeled $1$.
Similarly if the two loops are labeled $1$ and $X,$ then the merged loop must be labeled $X$ so that the local contribution for $\lambda$ goes from
$0$ to $-1.$ Finally, if the two loops are labeled $X$ and $X,$ then there is no label available for a single loop that will give $-3,$ so we define
$\partial$ to be zero in this case. We can summarize the result by saying that there is a multiplicative structure $m$ such that
$m(1,1) = 1, m(1,x) = m(x,1) = x, m(x,x) = 0,$ and this multiplication describes the structure of the partial differential when two loops merge.
Since this is the multiplicative structure of the algebra $V = k[x]/(x^{2}),$ we take this algebra as summarizing the differential.\\

Now consider the case where $s$ has a single loop at the site $\tau.$ Smoothing produces two loops. If the single loop is labeled $x,$ then we must label
each of the two loops by $x$ in order to make $\lambda$ decrease by $1$. If the single loop is labeled $1,$ then we can label the two loops by
$x$ and $1$ in either order. In this second case we take the partial differential of $s$ to be the sum of these two labeled states. This structure
can be described by taking a coproduct structure with $\Delta(x) = x \otimes x$ and 
$\Delta(1) = 1 \otimes x + x \otimes 1.$
We now have the algebra $V = k[x]/(x^{2})$ with product $m: V \otimes V \longrightarrow V$ and coproduct
$\Delta: V \longrightarrow V \otimes V,$ describing the differential completely. This completes the proof. //
\bigbreak

Partial differentials are defined on each enhanced state $s$ and a site $\tau$ of type
$A$ in that  state. We consider states obtained from the given state by  smoothing the given site $\tau$. The result of smoothing $\tau$ is to
produce a new state $s'$ with one more site of type $B$ than $s.$ Forming $s'$ from $s$ we either amalgamate two loops to a single loop at $\tau$, or
we divide a loop at $\tau$ into two distinct loops. In the case of amalgamation, the new state $s$ acquires the label on the amalgamated circle that
is the product of the labels on the two circles that are its ancestors in $s$. This case of the partial differential is described by the
multiplication in the algebra. If one circle becomes two circles, then we apply the coproduct. Thus if the circle is labeled $X$, then the resultant
two circles are each labeled $X$ corresponding to $\Delta(x) = x \otimes x$. If the orginal circle is labeled $1$ then we take the partial boundary
to be a sum of two enhanced states with  labels $1$ and $x$ in one case, and labels $x$ and $1$ in the other case,  on the respective circles. This
corresponds to $\Delta(1) = 1 \otimes x + x \otimes 1.$ Modulo two, the boundary of an enhanced state is the sum, over all sites of type $A$ in the
state, of the partial boundaries at these sites. It is not hard to verify directly that the square of the  boundary mapping is zero (this is the identity of
mixed partials!) and that it behaves
as advertised, keeping $j(s)$ constant. There is more to say about the nature of this construction with respect to Frobenius algebras and tangle
cobordisms. In Figures 4,5 and 6 we illustrate how the partial boundaries can be conceptualized in terms of surface cobordisms. Figure 4 shows how the partial boundary corresponds to a saddle point and illustrates the two cases of fusion and fission of circles. The equality of mixed
partials corresponds to topological equivalence of the corresponding surface cobordisms, and to the relationships between Frobenius algebras \cite{Kock} and the
surface cobordism category. In particular, in Figure 6 we show how in a key case of two sites (labeled 1 and 2 in that Figure) the two orders of partial
boundary are $$\partial_{2} \partial_{1} = (1 \otimes m) \circ (\Delta \otimes 1)$$ and
$$\partial_{1} \partial_{2} = \Delta \circ m.$$ In the Frobenius algebra $V = k[x]/(x^{2})$ we have the identity
$$(1 \otimes m) \circ (\Delta \otimes 1) = \Delta \circ m.$$ Thus the Frobenius algebra implies the identity of the mixed partials.
Furthermore, in Figure 5 we see that this identity corresponds to the topological equivalence of cobordisms under an exchange of saddle points.
\bigbreak

In Figures 7 and 8 we show another aspect of this algebra. As Figure 7 illustrates, we can consider
cup (minimum)  and cap (maximum)  cobordisms that go between the empty set and a single circle.
With the categorical arrow going down the page, the cap is a mapping from the base ring $k$ to 
the module $V$ and we denote this mapping by $\eta: k \longrightarrow V$. It is the {\em unit} for the algebra $V$ and is defined by $\eta(1) = 1_{V},$ taking $1$ in $k$ to $1_{V}$ in $V.$ The cup is a mapping from $V$ to $k$ and is denoted by $\epsilon: V \longrightarrow k.$ This is the {\em counit}.
As Figure 7 illustrates, we need a  basic identity about the counit which reads
$$\Sigma \epsilon(a_{1})a_{2} = a$$ for any $a \in V$ where 
$$\Delta(a) =  \Sigma a_{1} \otimes a_{2}.$$ The summation is over an appropriate set of elements in $v \otimes V$ as in our specific formulas for the algebra $k[x]/(x^{2}).$ Of course we also demand
$$\Sigma a_{1}\epsilon(a_{2}) = a$$ for any $a \in V.$ With these formulas about the counit and unit in
place, we see that cobordisms will give equivalent algebra when one cancels a maximum or a minimum with a saddle point, again as shown in Figure 7. 
\bigbreak

Note that for our algebra $V = k[x]/(x^{2}),$ it follows from the counit identies of the last paragraph that 
$$\epsilon(1) = 0$$ and $$\epsilon(x) = 1.$$ In fact, Figure 8 shows a formula that holds in this special algebra. The formula reads $$\epsilon(ab) = \epsilon(ax)\epsilon(b) + \epsilon(a)\epsilon(bx)$$
for any $a,b \in V.$ As the rest of Figure 8 shows, this identity means that a single tube in any cobordism can be cut, replacing it by a cups and a caps in a linear combination of two terms. The tube-cutting 
relation is shown in its most useful form at the bottom of Figure 8. In Figure 8, the black dots are 
symbols standing for the special element $x$ in the algebra.
\bigbreak

It is important to note that we have a nonsingular pairing $$\langle ~~|~~ \rangle: V \otimes V \longrightarrow k$$ defined by the equationn $$\langle a|b \rangle = \epsilon(ab).$$ One can define a 
Frobenius algebra by starting with the existence of a non-singular bilinear pairing.  In fact, a finite dimensional associative algebra with unit defined over a unital commutative ring  $k$ is said to be 
a {\em Frobenius algebra} if it is equipped with a non-degenerate bilinear form
$$\langle ~~|~~ \rangle: V \otimes V \longrightarrow k$$ such that
$$\langle ab | c \rangle = \langle a | bc \rangle$$ for all $a,b,c$ in the algebra. The other mappings and 
the interpretation in terms of cobordisms can all be constructed from this definition. See \cite{Kock}.
\bigbreak

\noindent{\bf Remark on Grading and Invariance.} In Section 2 we showed how the 
bracket, using the variable $q$, behaves under Reidemeister moves. These formulas correspond to how 
the invariance of the homology works in relation to the moves. We have that 
$$J_{K} = (-1)^{n_{-}} q^{n_{+} - 2n_{-}} \langle K \rangle,$$
where $n_{-}$ denotes the number of negative crossings in $K$ and $n_{+}$ denotes the number
of positive crossings in $K.$ $J(K)$ is invariant under all three Reidemeister moves. The corresponding formulas for Khonavov homology are as follows
$$J_{K} = (-1)^{n_{-}} q^{n_{+} - 2n_{-}} \langle K \rangle = $$
$$(-1)^{n_{-}} q^{n_{+} - 2n_{-}} \Sigma_{i,j} (-1)^{i}a^{j} dim(H^{i,j}(K) = $$
$$ \Sigma_{i,j} (-1) ^{i + n_{+} }q^{j + n_{+} - 2n_{-1} }dim(H^{i,j}(K)) = $$
$$\Sigma_{i,j} (-1)^{i} q^{j} dim(H^{i - n_{-}, j - n_{+} + 2n_{-}}(K)).$$
It is often more convenient to define the {\em Poincar\'e polynomial} for Khovanov homology via
$$P_{K}(t, q)  
= \Sigma_{i,j} t^{i} q^{j} dim(H^{i - n_{-}, j - n_{+} + 2n_{-}}(K)).$$ 
The Poincar\'e polynomial is a two-variable polynomial invariant of knots and links, generalizing the Jones polynomial. Each coefficient $$ dim(H^{i - n_{-}, j - n_{+} + 2n_{-}}(K))$$ is an 
invariant of the knot, invariant under all three Reidemeister moves. In fact, the homology groups 
$$H^{i - n_{-}, j - n_{+} + 2n_{-}}(K)$$ are knot invariants. The grading compensations show how the grading of the homology can change from diagram to diagram for diagrams that represent the same knot.
\bigbreak

\begin{figure}
     \begin{center}
     \begin{tabular}{c}
     \includegraphics[width=6cm]{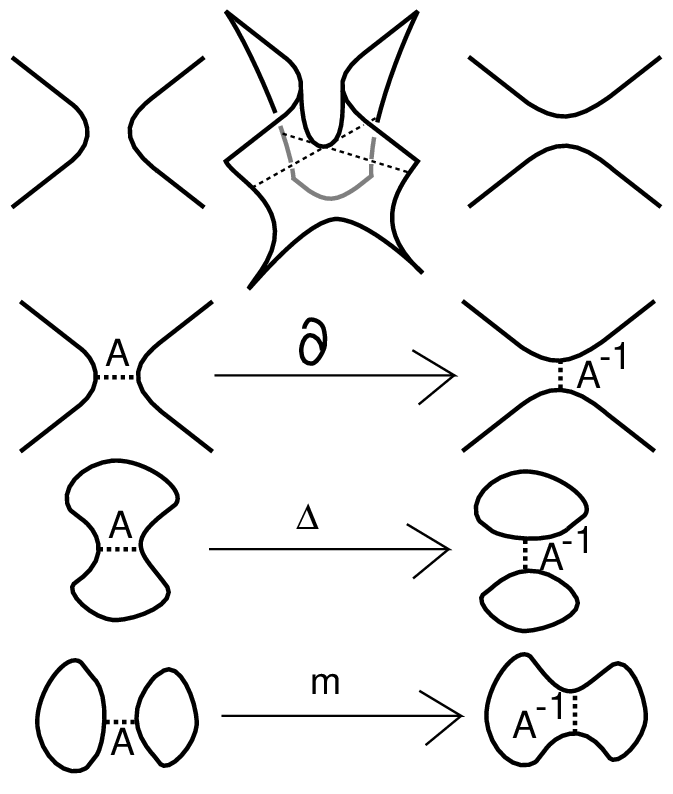}
     \end{tabular}
     \caption{\bf SaddlePoints and State Smoothings}
     \label{Figure 4}
\end{center}
\end{figure}

\begin{figure}
     \begin{center}
     \begin{tabular}{c}
     \includegraphics[width=8cm]{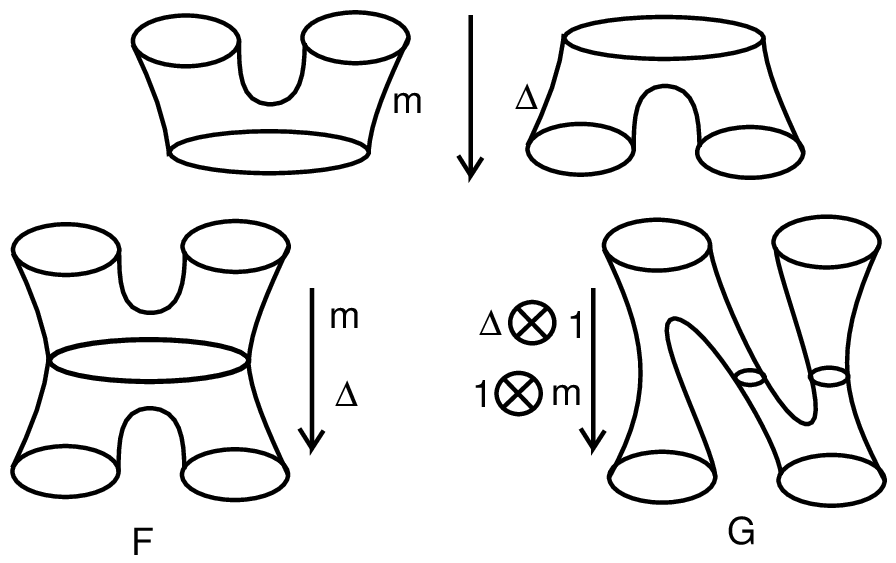}
     \end{tabular}
     \caption{\bf Surface Cobordisms}
     \label{Figure 5}
\end{center}
\end{figure}

\begin{figure}
     \begin{center}
     \begin{tabular}{c}
     \includegraphics[width=9cm]{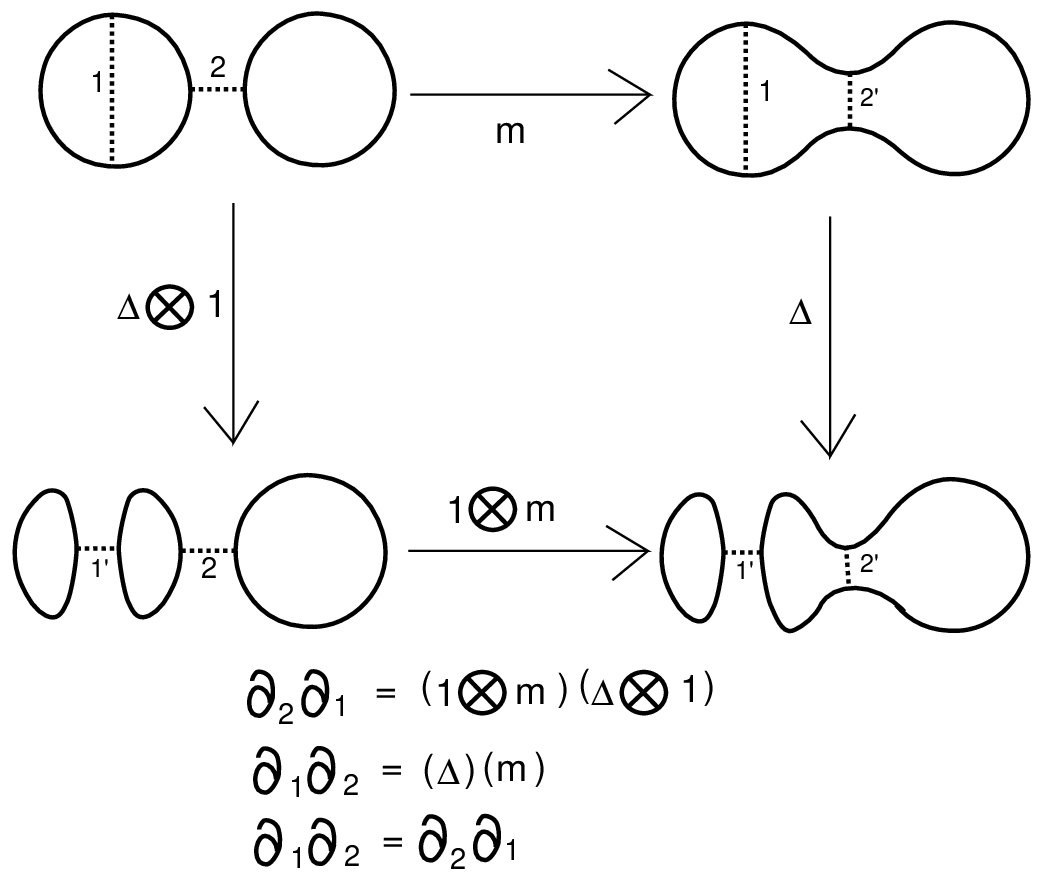}
     \end{tabular}
     \caption{\bf Local Boundaries Commute}
     \label{Figure 6}
\end{center}
\end{figure}

\begin{figure}
     \begin{center}
     \begin{tabular}{c}
     \includegraphics[width=8cm]{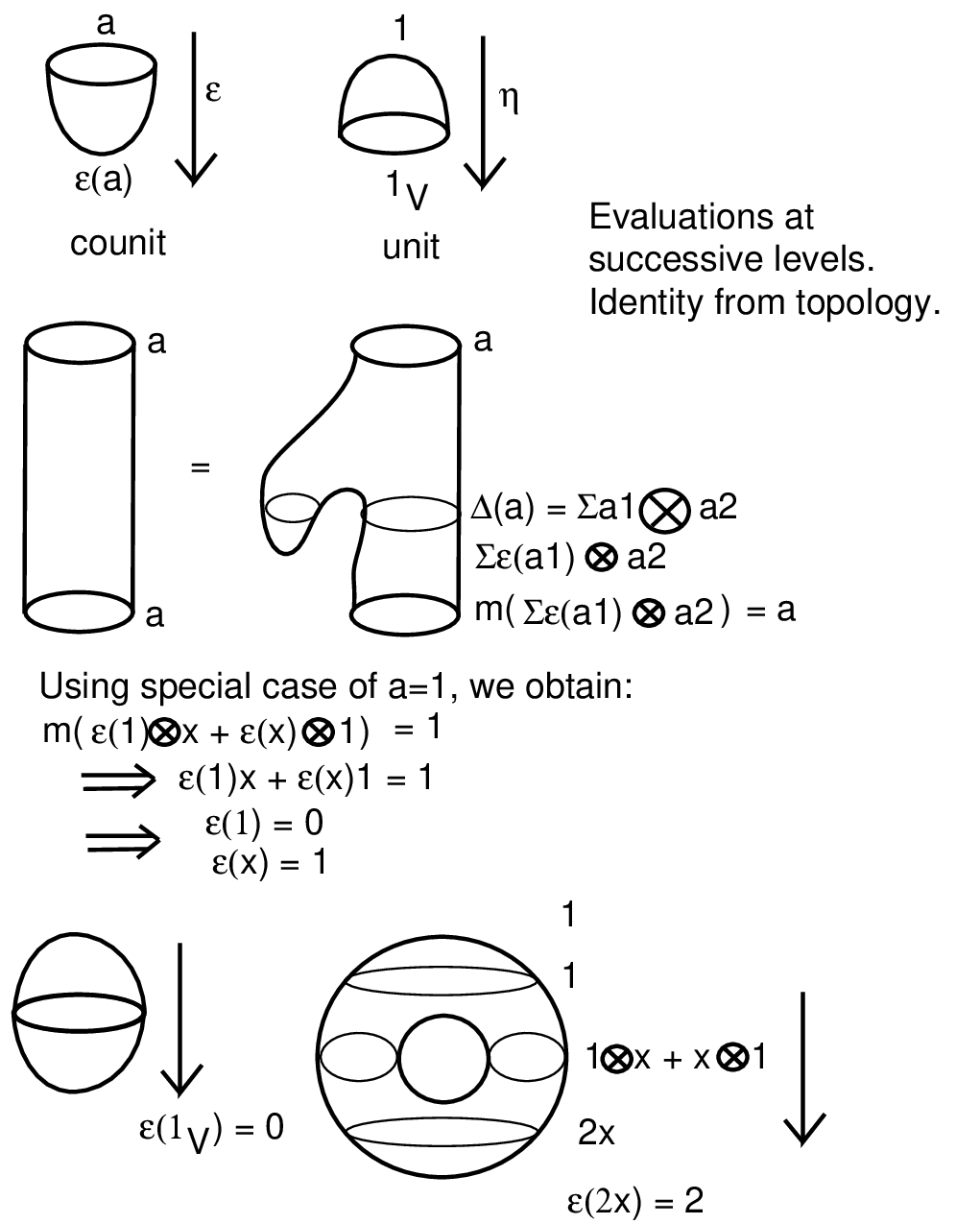}
     \end{tabular}
     \caption{\bf Unit and Counit as Cobordisms}
     \label{Figure 7}
\end{center}
\end{figure}

\begin{figure}
     \begin{center}
     \begin{tabular}{c}
     \includegraphics[width=8cm]{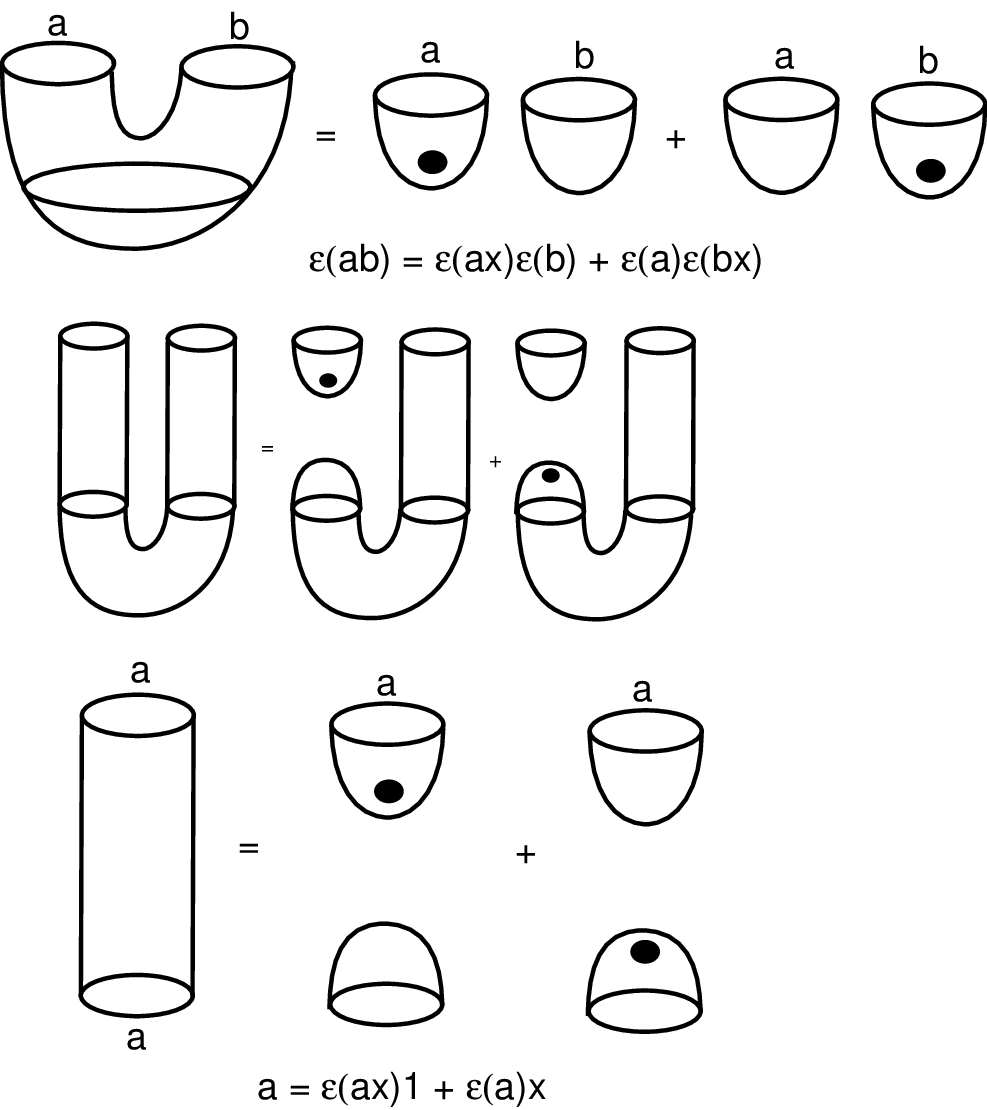}
     \end{tabular}
     \caption{\bf The Tube Cutting Relation}
     \label{Figure 8}
\end{center}
\end{figure}

\noindent {\bf Remark on Integral Differentials}. Choose an ordering for the crossings in the link diagram $K$ and denote them by $1,2,\cdots n.$
Let $s$ be any enhanced state of $K$ and let $\partial_{i}(s)$ denote the chain obtained from $s$ by applying a partial boundary at the $i$-th site
of $s.$ If the $i$-th site is a smoothing of type $A^{-1}$, then $\partial_{i}(s) = 0.$ If the $i$-th site is a smoothing of type $A$, then
$\partial_{i}(s)$ is given by the rules discussed above (with the same signs). The compatibility conditions that we have discussed show that
partials commute in the sense that $\partial_{i} (\partial_{j} (s)) = \partial_{j} (\partial_{i} (s))$ for all $i$ and $j.$ One then defines
signed boundary formulas in the usual way of algebraic topology. One way to think of this regards the complex as the analogue of
a complex in de Rham cohomology. Let $\{dx_{1}, dx_{2},\cdots, dx_{n}\}$ be a formal basis for a Grassmann algebra so that
$dx_{i} \wedge dx_{j} = - dx_{j} \wedge dx_{i}$ Starting with enhanced states $s$ in $C^{0}(K)$ (that is, states with all $A$-type smoothings)
define formally, $d_{i}(s) = \partial_{i}(s) dx_{i}$ and regard $d_{i}(s)$ as identical with $\partial_{i} (s)$ as we have previously regarded it
in $C^{1}(K).$ In general, given an enhanced state $s$ in $C^{k}(K)$ with $B$-smoothings at locations $i_{1} < i_{2} < \cdots < i_{k},$ we represent
this chain as $s \, dx_{i_{1}} \wedge \cdots \wedge dx_{i_{k}}$ and define
$$\partial ( s \, dx_{i_{1}} \wedge \cdots \wedge dx_{i_{k}} ) = \sum_{j=1}^{n} \partial_{j}(s) \, dx_{j} \wedge dx_{i_{1}} \wedge \cdots \wedge dx_{i_{k}},$$
just as in a de Rham complex. The Grassmann algebra automatically computes the correct signs in the chain complex, and this boundary formula gives the
original boundary formula when we take coefficients modulo two. Note, that in this formalism, partial differentials $\partial_{i}$ of enhanced states with a
$B$-smoothing at the site $i$ are zero due to the fact that $dx_{i} \wedge dx_{i} = 0$ in the Grassmann algebra. There is more to discuss
about the use of Grassmann algebra in this context. For example, this approach clarifies parts of the construction in \cite{M1}.
\bigbreak

It of interest to examine this analogy between the Khovanov (co)homology and de Rham cohomology. In that analogy the enhanced
states correspond to the differentiable functions on a manifold. The Khovanov complex $C^{k}(K)$ is generated by elements of the form
$s \, dx_{i_{1}} \wedge \cdots \wedge dx_{i_{k}}$ where the enhanced state $s$ has $B$-smoothings at exactly the sites $i_{1},\cdots, i_{k}.$
If we were to follow the analogy with de Rham cohomology literally, we would define a new complex $DR(K)$ where $DR^{k}(K)$ is generated by elements
$s \, dx_{i_{1}} \wedge \cdots \wedge dx_{i_{k}}$ where $s$ is {\it any} enhanced state of the link $K.$ The partial boundaries are defined in the same
way as before and the global boundary formula is just as we have written it above. This gives a {\it new} chain complex associated with the link $K.$
Whether its homology contains new topological information about the link $K$ will be the subject of a subsequent paper.
\bigbreak

In the case of de Rham cohomology, we can also look for compatible unitary transformations. Let $M$ be a differentiable manifold and ${\mathcal
C}(M)$ denote the DeRham complex of
$M$ over the complex numbers. Then for a differential form of the type $f(x) \omega$ in local coordinates $x_{1}, \cdots ,x_{n}$ and $\omega$ a wedge
product of a subset of
$dx_{1} \cdots dx_{n},$ we have
$$d(f \omega) = \sum_{i=1}^{n} (\partial f/\partial x_{i}) dx_{i} \wedge \omega.$$ Here $d$ is the differential for the DeRham complex. Then
${\mathcal C}(M)$ has as basis the set of $|f(x) \omega \rangle$ where $\omega = dx_{i_{1}} \wedge \cdots \wedge dx_{i_{k}}$ with
$i_{1} < \cdots < i_{k} .$ We could achieve $U d + d U = 0$ if $U$ is a very simple unitary operator (e.g. multiplication by phases that do not depend on
the coordinates $x_{i}$) but in general it will be an interesting problem to determine all unitary operators $U$ with this property.
\bigbreak

\noindent {\bf A further remark on de Rham cohomology.} There is another relation with the de
Rham complex: In \cite{Pr} it was observed that Khovanov homology is related to Hochschild
homology and Hochschild homology is thought to be an algebraic version of de Rham chain
complex (cyclic cohomology corresponds to de Rham cohomology), compare \cite{Lo}.
\bigbreak

\section {\bf The Cube Category and the Tangle Cobordism Structure of Khovanov Homology}
We can now connect the constructions of the last  section
with the homology construction via the cube category. Here it will be convenient to think of the state category ${\mathcal S}(K)$ as a cube category with extra structure. Thus we will denote the 
bracket states by sequences of $A$'s and $B$'s as in Figures 2 and 3. And we shall regard the
maps such as $d_{2}: \langle AABA \rangle \longrightarrow \langle ABBA \rangle$ as corresponding
to re-smoothings of bracket states that either join or separate state loops. 
We take $V = k[x]/(x^{2})$ with the coproduct
structure as given in the previous section. The maps from $m: V \otimes V \longrightarrow V$ and 
$\Delta: V \longrightarrow V \otimes V$ allow us to define the images of the resmoothing maps 
$d_{i}$ under a functor
${\mathcal F}: {\mathcal S}(K) \longrightarrow {\mathcal M}$ where ${\mathcal M}$ is the category generated by $V$ by taking  tensor powers of $V$  and direct sums of these tensor powers. It then follows that the homology we have described in the previous section is exactly the homology associated with this functor.
\bigbreak

The material in the previous section also suggests a modification of the state category  
${\mathcal S}(K).$ Instead of taking the maps in this category to be simply the abstract arrows generated by elementary re-smoothings of states from $A$ to $B,$ we can regard each such smoothing as a 
surface cobordism from the set of circles comprising the domain state to the set of circles comprising the codomain state. With this, in mind, two such cobordisms represent equivalent morphisms whenever
the corresponding surfaces are homeomorphic relative to their boundaries. Call this category
${\mathcal CobS}(K).$ We then easily generalize the observations of the previous section, particularly
Figures 4, 5 and 6, to see that we have the desired commuting relations $d_{i} d_{j} = d_{j}d_{i}$
(for $i \ne j$) in  ${\mathcal CobS}(K)$ so that any functor from ${\mathcal CobS}(K)$ to a module 
category will have a well-defined chain complex and associated homology. This applies, in particular to 
the functor we have constructed, using the Frobenius algebra $V = k[x]/(x^{2}).$
\bigbreak

In \cite{BN} BarNatan takes the approach using surface cobordisms a step further by making a categorical analog of the chain complex. For this purpose we let ${\mathcal CobS}(K)$ become an 
additive category. Maps between specific objects $X$ and $Y$ added formally and the set
$Maps(X,Y)$ is a module over the integers. More generally, let ${\mathcal C}$ be an additive category.
In order to create the analog of a chain complex, let ${\mathcal Mat(C)}$ denote the {\em Matrix Category of $C$} whose objects are $n$-tuples (vectors) of objects of ${\mathcal C}$  ($n$ can be any natural number) and whose morphisms are in the form of a  matrix  $m = (m_{ij})$ of morphisms in $C$ where we write
$$m : O \longrightarrow O'$$
and $$m_{ij}: O_{i} \longrightarrow O'_{j}$$
for 
$$O = (O_{1},\cdots,O_{n}),$$ 
$$O' = (O'_{1},\cdots,O'_{m}).$$
Here $O_{i}$ and $O'_{j}$ are objects in ${\mathcal C}$ while $O$ and $O'$ are objects in 
${\mathcal Mat(C)}.$ Composition of morphisms in ${\mathcal Mat(C)}$ follows  the pattern of 
matrix multiplication. If   $$n : O' \longrightarrow O''$$ then
$$n \circ m : O \longrightarrow O'' $$ and $$(n \circ m)_{i,j} = \Sigma_{k} n_{i,k} \circ m_{k,j}$$
where the compositions in the summation occur in the category $C.$
\bigbreak

We then define the {\em category of complexes over ${\mathcal C}$}, denoted
${\mathcal Kom(Mat(C))}$ to consist of sequences of objects of ${\mathcal Mat(C)}$ and maps between
them so that consecutively composed maps are equal to zero. 
$$ \cdots \longrightarrow O^{k} \longrightarrow O^{k+1}  \longrightarrow
 O^{k+2} \longrightarrow \cdots .$$ 
 Here we let $\partial_{k}:O^{k} \longrightarrow O'^{k+1}$
 denote the differential in the complex and we assme that $\partial_{k+1} \partial_{k} = 0.$
 A morphism between complexes $O^{*}$ and $O'^{*}$ consists in a family of maps
 $f_{k} : O^{k} \longrightarrow O'^{k}$ such that $\partial'_{k} f_{k} = f_{k+1} \partial_{k}.$  
 Such morphisms will be called {\em chain maps}.
 \bigbreak

At this abstract level, we cannot calculate homology since kernels and cokernels are not available, but we can define the {\em homotopy type} of a complex in ${\mathcal Kom(Mat(C))}.$
We say that two chain maps $f:O \longrightarrow O'$ and $g:O \longrightarrow O'$ are {\em homotopic} if there is a sequence of mappings $H_{k}:O^{k} \longrightarrow O^{k-1}$ such that 
$$f - g = H \partial + \partial H.$$ Specifically, this means that 
$$f_{k} - g_{k} = H_{k+1} \partial_{k} + \partial_{k+1} H_{k}.$$
Note that if $\phi = H \partial + \partial H,$ then $$\partial \phi = \partial H \partial = \phi \partial.$$
Thus any such $\phi$ is a chain map. We call two complexes $O$ and $O'$ {\em homotopy equivalent}
if there are chain maps $F: O \longrightarrow O'$ and $G:O' \longrightarrow O$ such that both
$FG$ and $GF$ are homotopic to the identity map of $O$ and $O'$ respectively. The homotopy type of a complex is an abstract substitute for the homology since, in an abelian category (where one can compute homology) homology is an invariant of homotopy type.
\bigbreak

We are now in a position to work with the category ${\mathcal Kom(Mat(CobS(K)))}$ where $K$ is a 
link diagram. The question is, what extra equivalence relation on the category 
 ${\mathcal CobS(K)}$ will ensure that the homotopy types in ${\mathcal Kom(Mat(CobS(K)))}$
 will be invariant under Reidemeister moves on the diagram $K.$ 
 \bigbreak
 
 BarNatan \cite{BN} gives an elegant answer to this question. His answer is illustrated in Figure 9 where we show the {\em 4Tu Relation}, the {\em Sphere Relation} and the {\em Torus Relation}. 
 {\it The key relation is the the $4Tu$ relation.} The  $4Tu$ relation
 serves a number of purposes, including being a basic homotopy in the category
  ${\mathcal Kom(Mat(CobS(K)))}.$ \\
  
  The $4Tu$ relation can be described as follows: There are 
  four local bits of surface, call them $S_{1},S_{2},S_{3},S_{4}.$ Let $C_{i,j}$ denote this configuration 
  with a tube connecting $S_{i}$ and $S_{j}.$ Then in the cobordism category we take the identity
  $$C_{1,2} + C_{3,4} = C_{1,3} + C_{2,4}.$$ It is a good exercise for the reader to show that the 
  $4Tu$ relation follows from the tube cutting relation of Figure 8. In fact Figures 15 gives a schematic for the four-term relation, where arrows correspond to tubes attached to surfaces and arcs correspond to surfaces. \\
  
  Figure 16 shows how the tube-cutting relation is a  consequence of the  $4Tu$ relation, when it is assumed that the chain homotopy theory occurs over a ring where $2$ is an invertible element.
  Without this assumption, we cannot perform the trick, indicated in Figure 16, of packing up a punctured torus (divided by $2$) as a ``dot".  This dot will later be interpreted (in the next section) as an element in an algebra. If $2$ is not invertible then there is no translation of the $4Tu$ relation to at tube--cutting relation and the chain- homotopy theory will be different. {\it For the remainder of this paper, we assume that $2$ is invertible.} Figure 17 shows a derivation of the  $4Tu$ relation from the tube cutting relation.\\
  
Note that the Sphere and Torus relations assert  that the $2$--sphere has value $0$ and that the torus has value $2,$ just as we have seen by using the Frobenius algebra in Figure 7.\\

 To illustrate how things work once we factor by these relations, we show in Figures 10 and 11 how one sees parts the homotopy equivalance of the complexes for a diagram before and after the second Reidemeister move. In Figure 10 we show the complexes and indicate chain maps $F$ and $G$ between them and homotopies in the complex for the diagram before it is simplified by the Reidemeister move. In Figures 11 and 12 we show some of the cobordism compositions of the maps in this complex. In  Figure 13 we show these maps and their compositions in the form of a four-term identity that 
 verifies the needed chain-homotopy for the equivalence of the complexes before and after the Reidemeister move. Figure 14 shows the same pattern as Figure 13, but is designed to make it clear that this identity is indeed exactly the $4Tu$ relation!
 { \it Thus the $4Tu$ relation is the key to the chain-homotopy invariance of the Khovanov Complex under the Second Reidemeister move.}\\
 
 As shown in Figure 13, each of the terms in the relation is factored into 
 mappings involving $F_1$, $G_1$ and the homotopies $H_{1}$ and  $H_{2}$ and the boundary
 mappings in the complex. Study of Figure 13 will convince the reader that the complexes before
 and after the second Reidemeister move are homotopy equivalent. A number of details are left to the reader. For example, note that in Figure 10 we have indicated the categorical chain complexes
 $Z$ and $W$ by showing only how they differ locally near the change corresponding to a Reidemeister two move. We give, via Figures 10 and 11, chain maps $F:W \longrightarrow Z$ and 
 $G:Z  \longrightarrow W.$ These maps consist in a particular cobordism on one part of the complex and an identity map on the other part of the complex. We have specifically labeled parts of these mappings by $F_1$ and by $G_1.$ Using the implicit definitions of $F_1$ and $G_1$ given in Figure 11, the reader will easily see that 
 $G_{1}F_{1} = 0$ since this composition includes a $2$-sphere. From this it follows that $GF$ is the identity mapping on the complex $W.$ We also leave to the reader to check that the mappings $F$ and $G$  commute with the boundary mappings so that they are mappings of complexes. The part of the homotopy indicated shows that $FG$ is homotopic to the identity (up to sign) and so shows that the complexes $Z$ and $W$ are homotopy equivalent. One needs the value of the torus equal to $2$ for homotopy invariance under the first Reidemeister move.
 Invariance under the third Reidemeister move can be deduced from invariance under the second Reidemeister move and a description of the (abstract) chain complex ${\cal C}(\Across)$
 as the mapping cone of ${\cal C}(\Asmooth) \longrightarrow {\cal C}(\Bsmooth)$ in a direct generalization of the original argument that shows that the bracket polynomial is invariant under the third Reidemeister move
 as a consequence of its invariance under the second Reidemeister move.
This is the main part of the  full derivation of homotopy equivalences corresponding to all three Reidemeister moves that is given in \cite{BN}. 
 \bigbreak

\begin{figure}
     \begin{center}
     \begin{tabular}{c}
     \includegraphics[width=9cm]{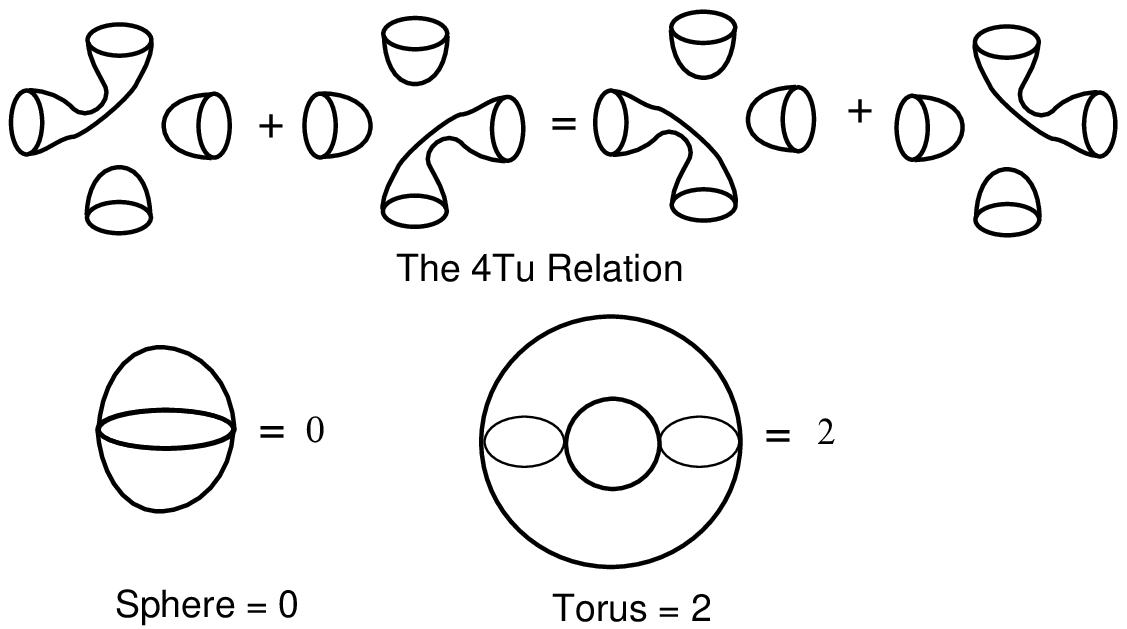}
     \end{tabular}
     \caption{\bf The 4Tu Relation, Sphere and Torus Relations}
     \label{Figure 9}
\end{center}
\end{figure}

\begin{figure}
     \begin{center}
     \begin{tabular}{c}
     \includegraphics[width=9cm]{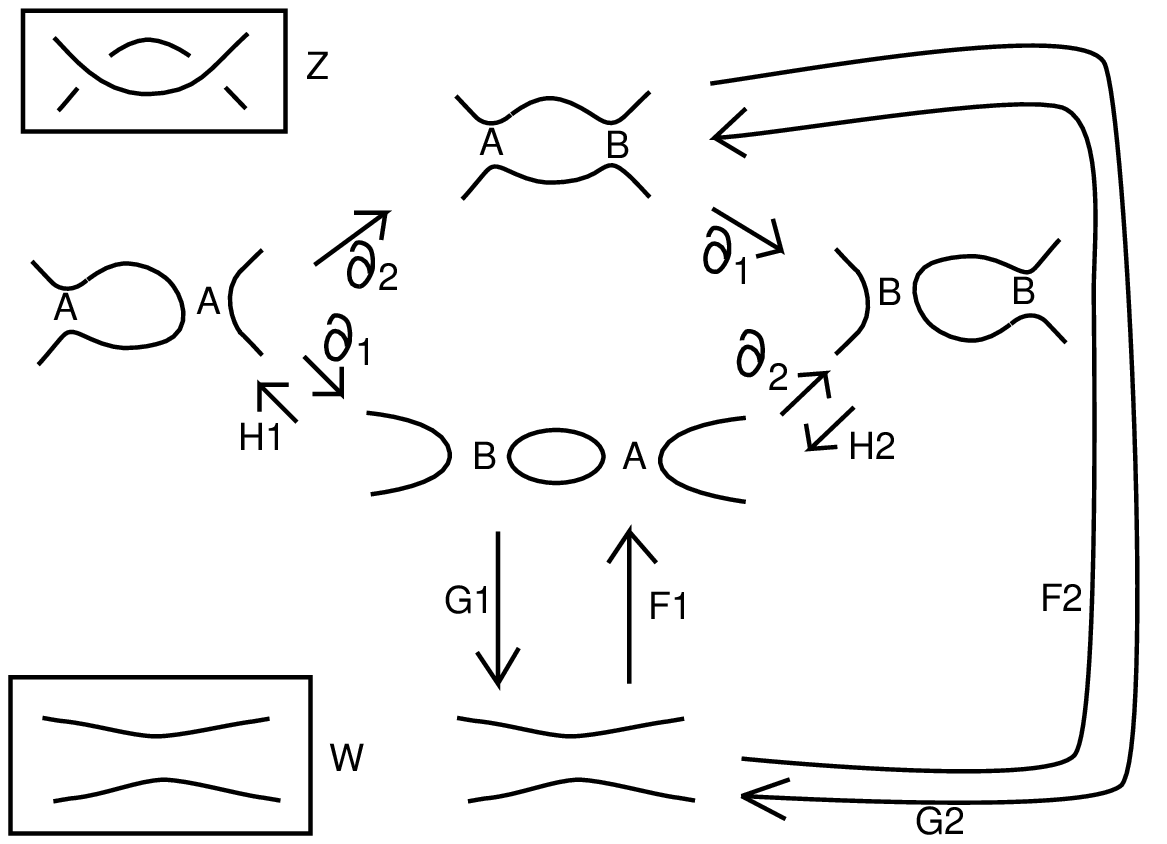}
     \end{tabular}
     \caption{\bf Complexes for Second Reidemeister Move}
     \label{Figure 10}
\end{center}
\end{figure}

\clearpage

\begin{figure}
     \begin{center}
     \begin{tabular}{c}
     \includegraphics[width=9cm]{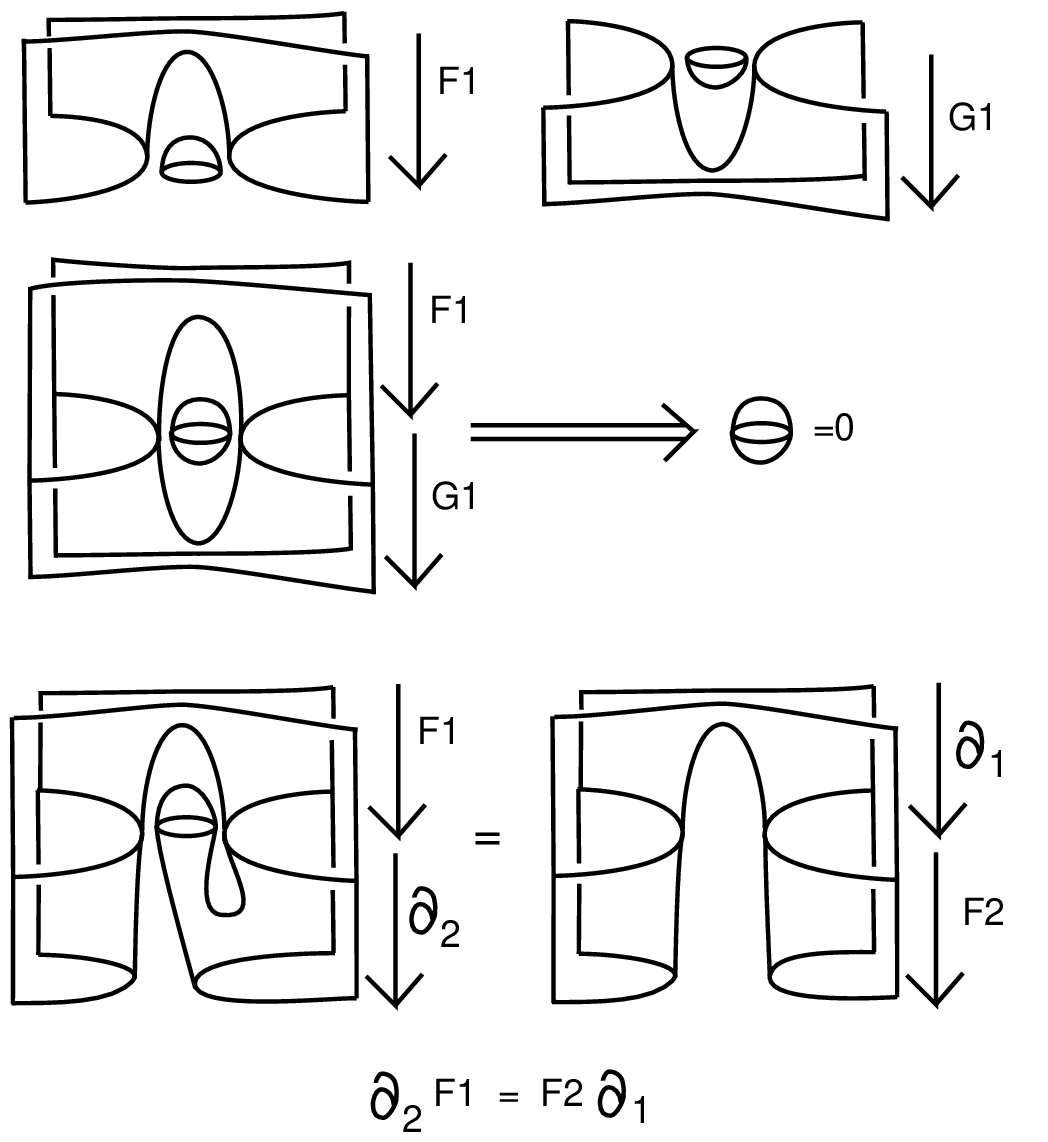}
     \end{tabular}
     \caption{\bf Cobordism Compositions for Second Reidemeister Move}
     \label{cobmaps}
\end{center}
\end{figure}

\begin{figure}
     \begin{center}
     \begin{tabular}{c}
     \includegraphics[width=9cm]{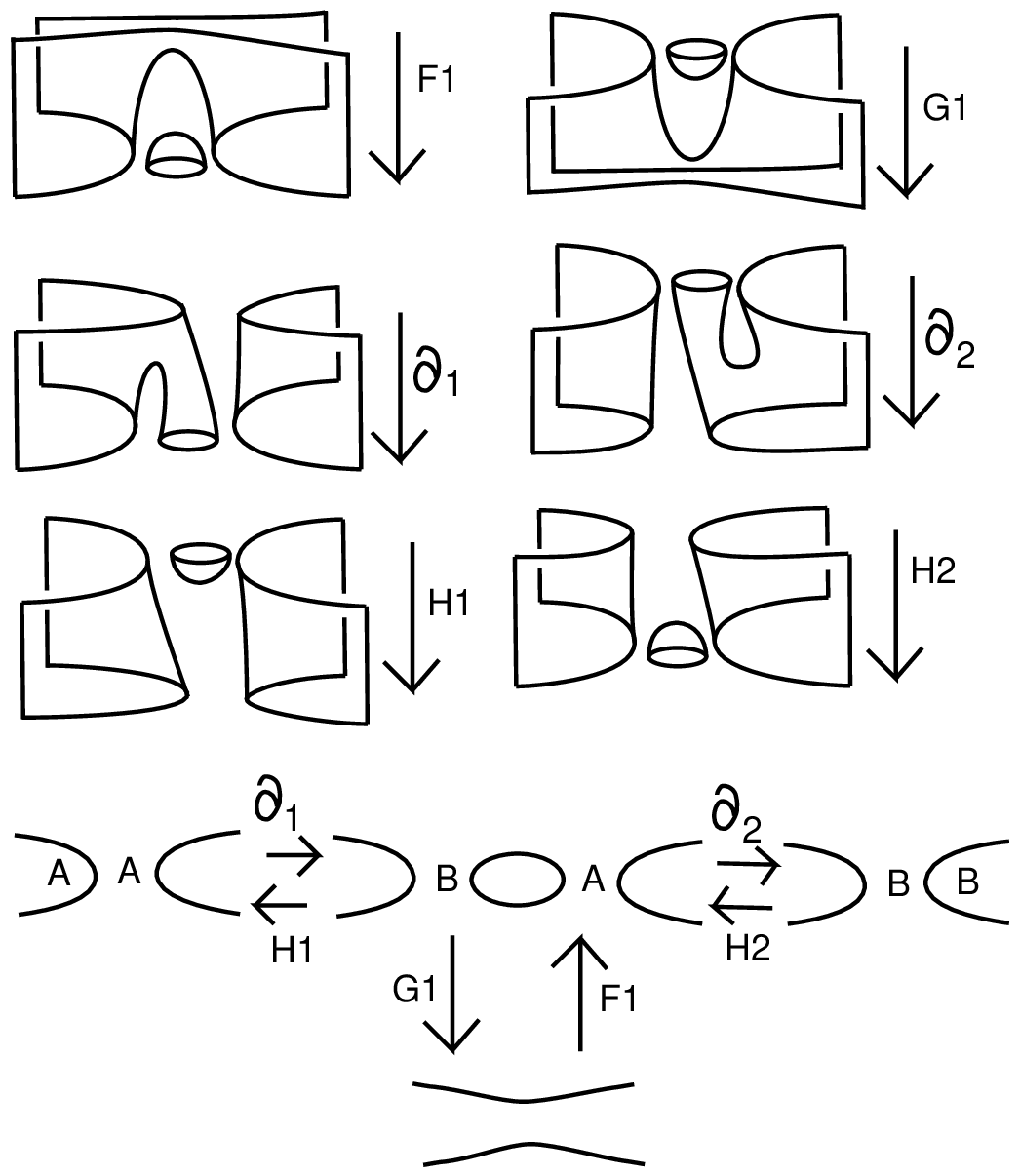}
     \end{tabular}
     \caption{\bf Preparation for Homotopy for Second Reidemeister Move}
     \label{cobmaps1}
\end{center}
\end{figure}

\begin{figure}
     \begin{center}
     \begin{tabular}{c}
     \includegraphics[width=9cm]{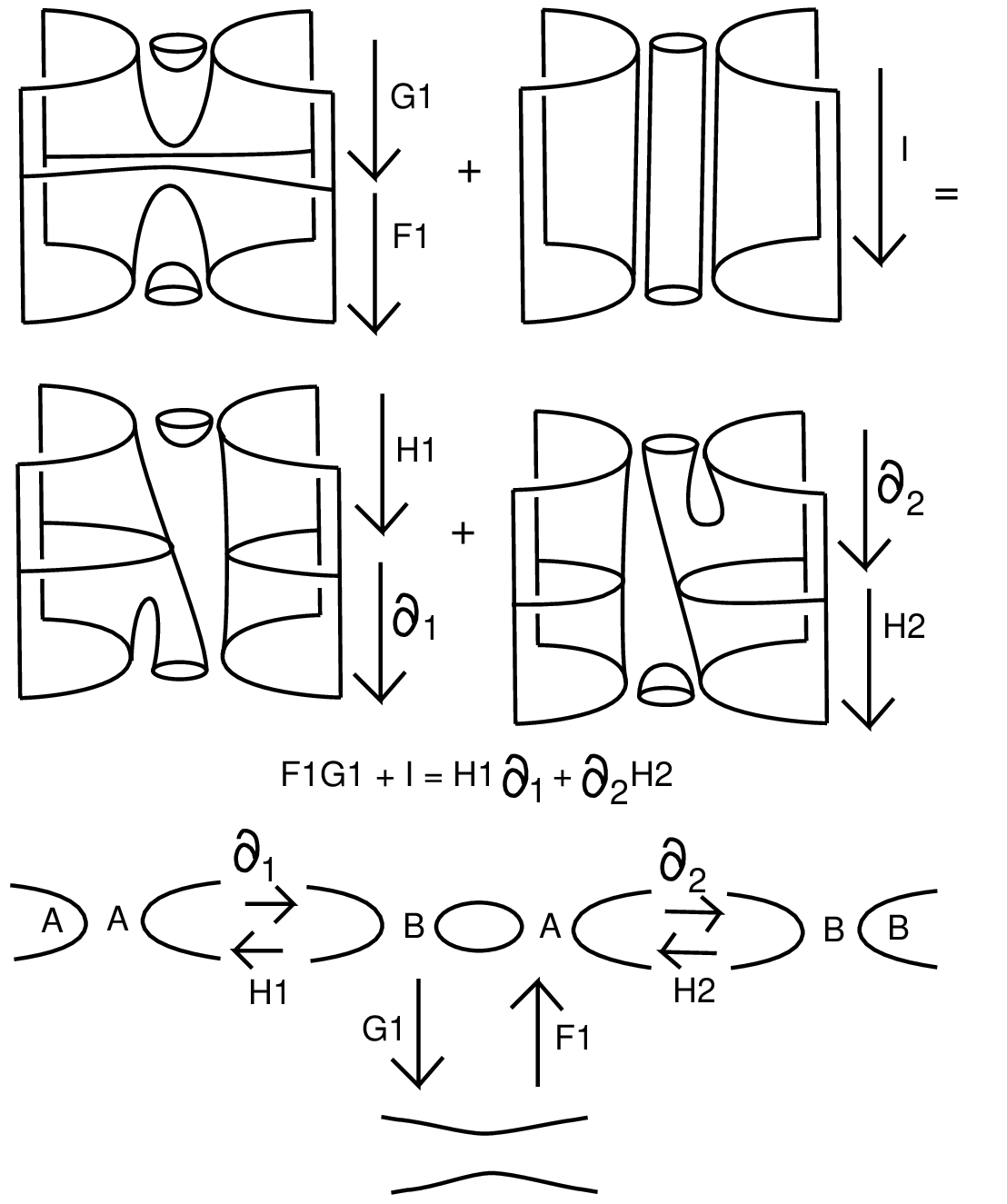}
     \end{tabular}
     \caption{\bf Homotopy for Second Reidemeister Move}
     \label{Figure 11}
\end{center}
\end{figure}

\begin{figure}
     \begin{center}
     \begin{tabular}{c}
     \includegraphics[width=9cm]{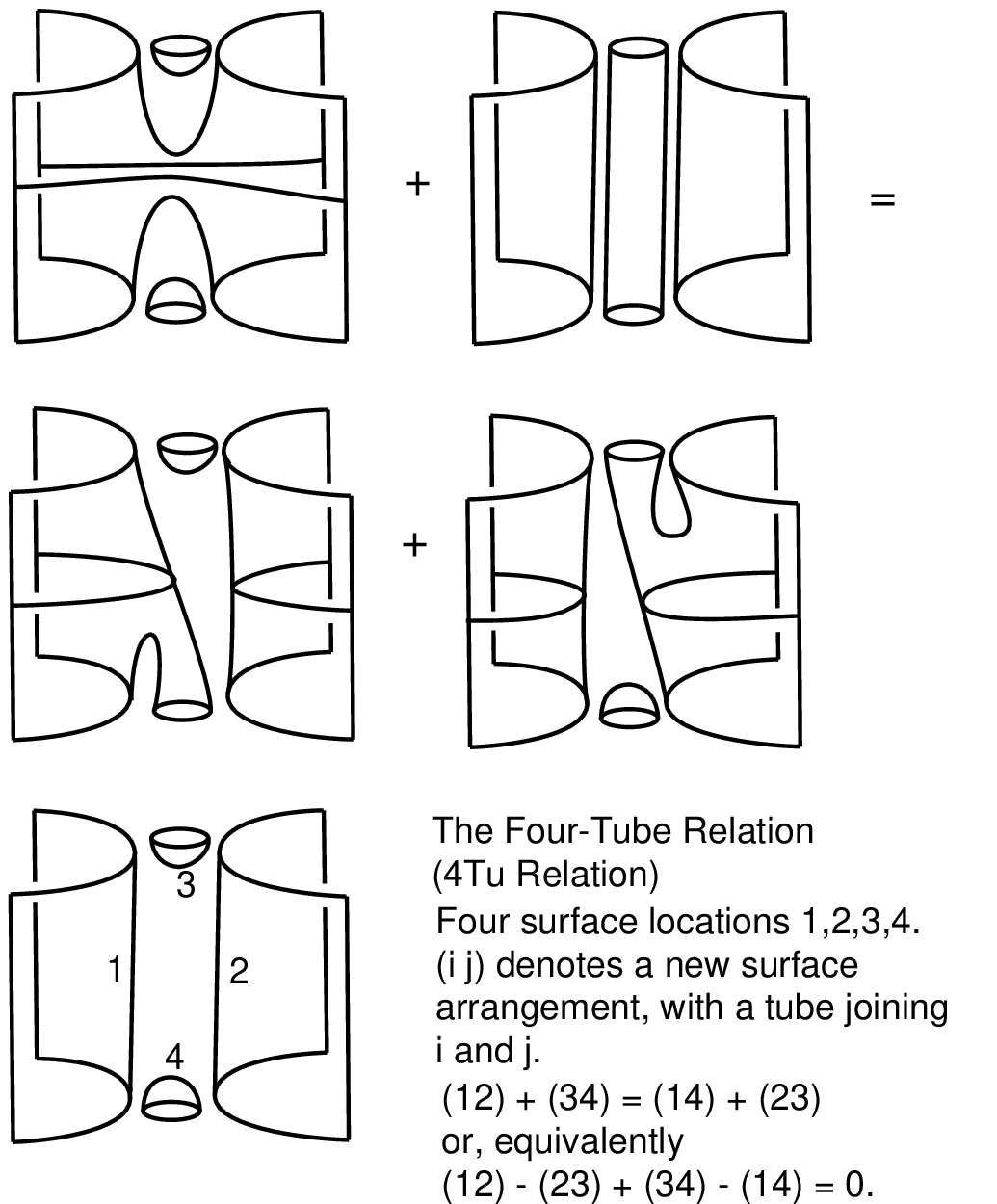}
     \end{tabular}
     \caption{\bf Four-Tube  Relation From Homotopy}
     \label{surfacerelation}
\end{center}
\end{figure}

\begin{figure}
     \begin{center}
     \begin{tabular}{c}
     \includegraphics[width=9cm]{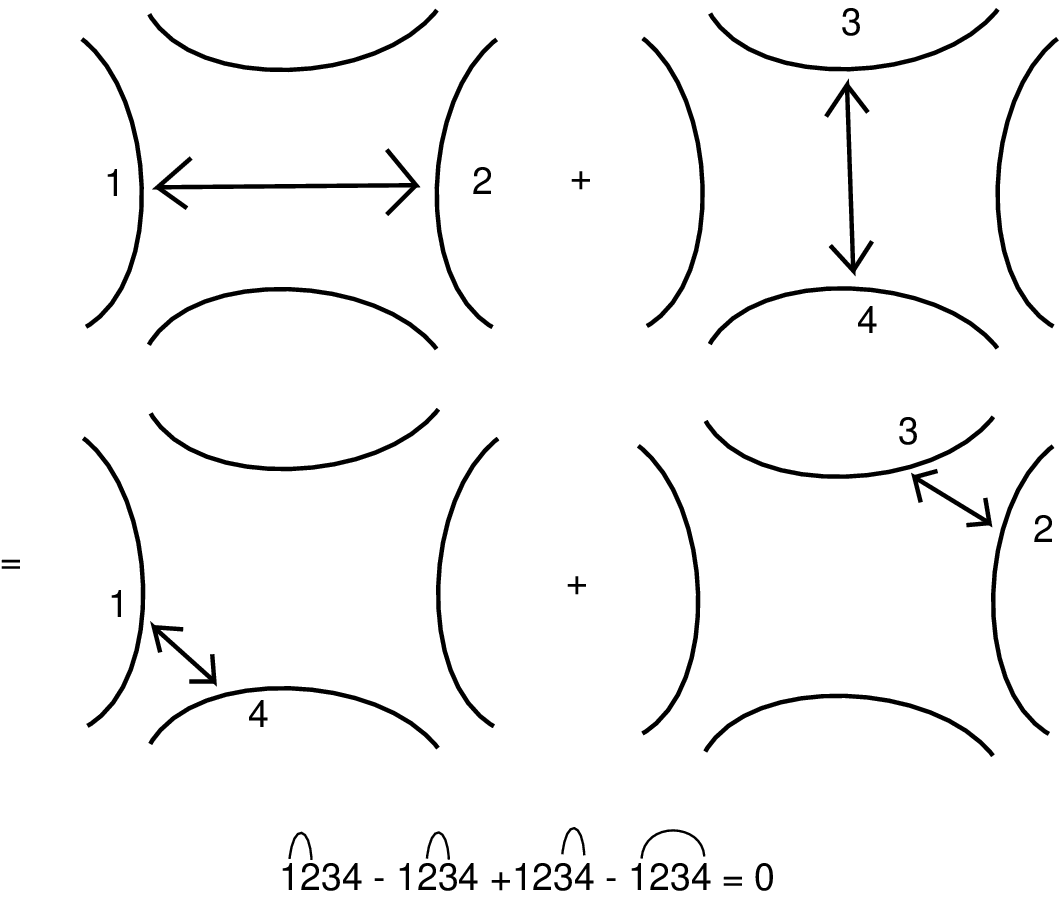}
     \end{tabular}
     \caption{\bf Schematic Four-Tube Relation}
     \label{schematic}
\end{center}
\end{figure}

\begin{figure}
     \begin{center}
     \begin{tabular}{c}
     \includegraphics[width=9cm]{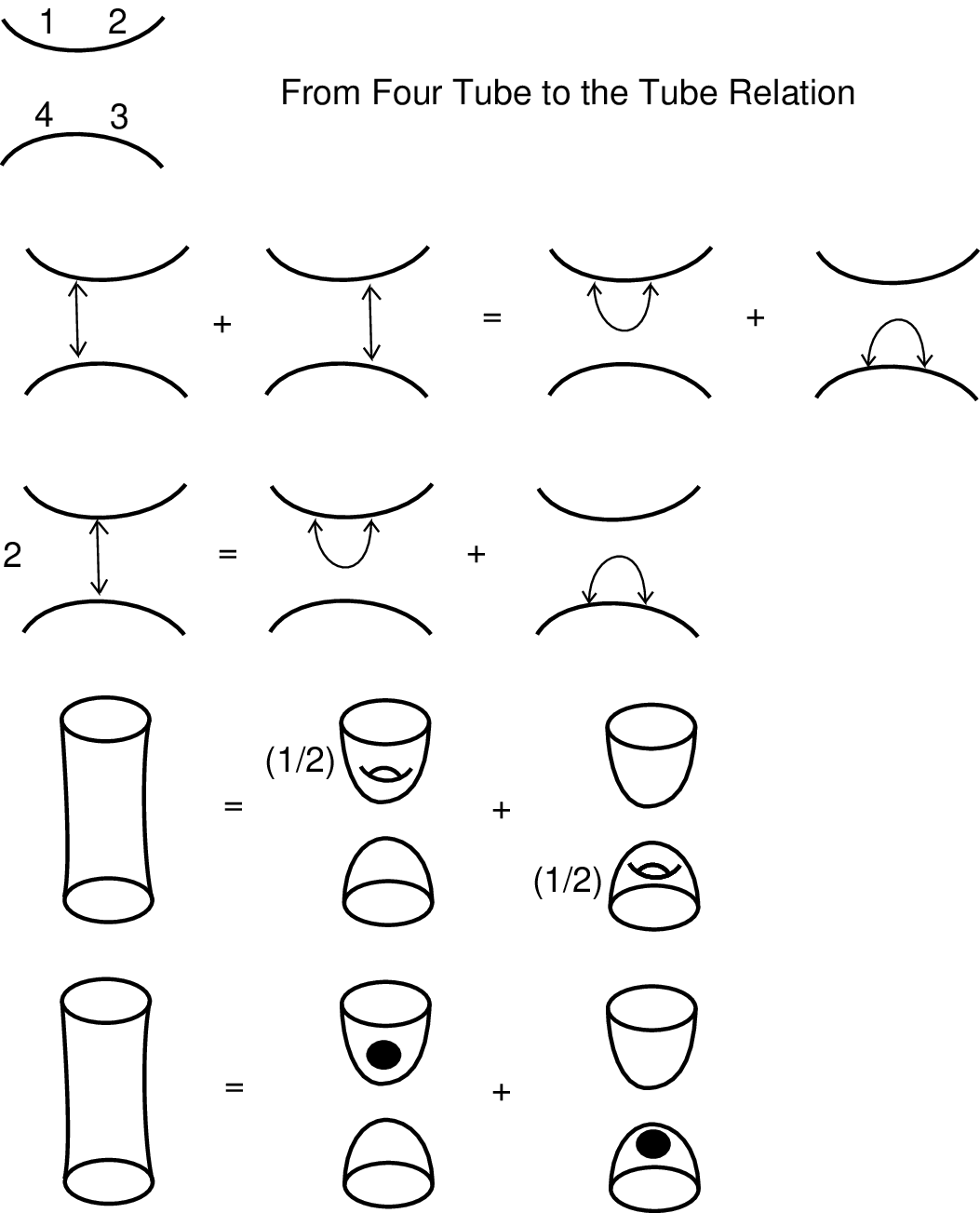}
     \end{tabular}
     \caption{\bf From Four-Tube Relation to Tube-Cutting Relation}
     \label{fourtotu}
\end{center}
\end{figure}

\begin{figure}
     \begin{center}
     \begin{tabular}{c}
     \includegraphics[width=9cm]{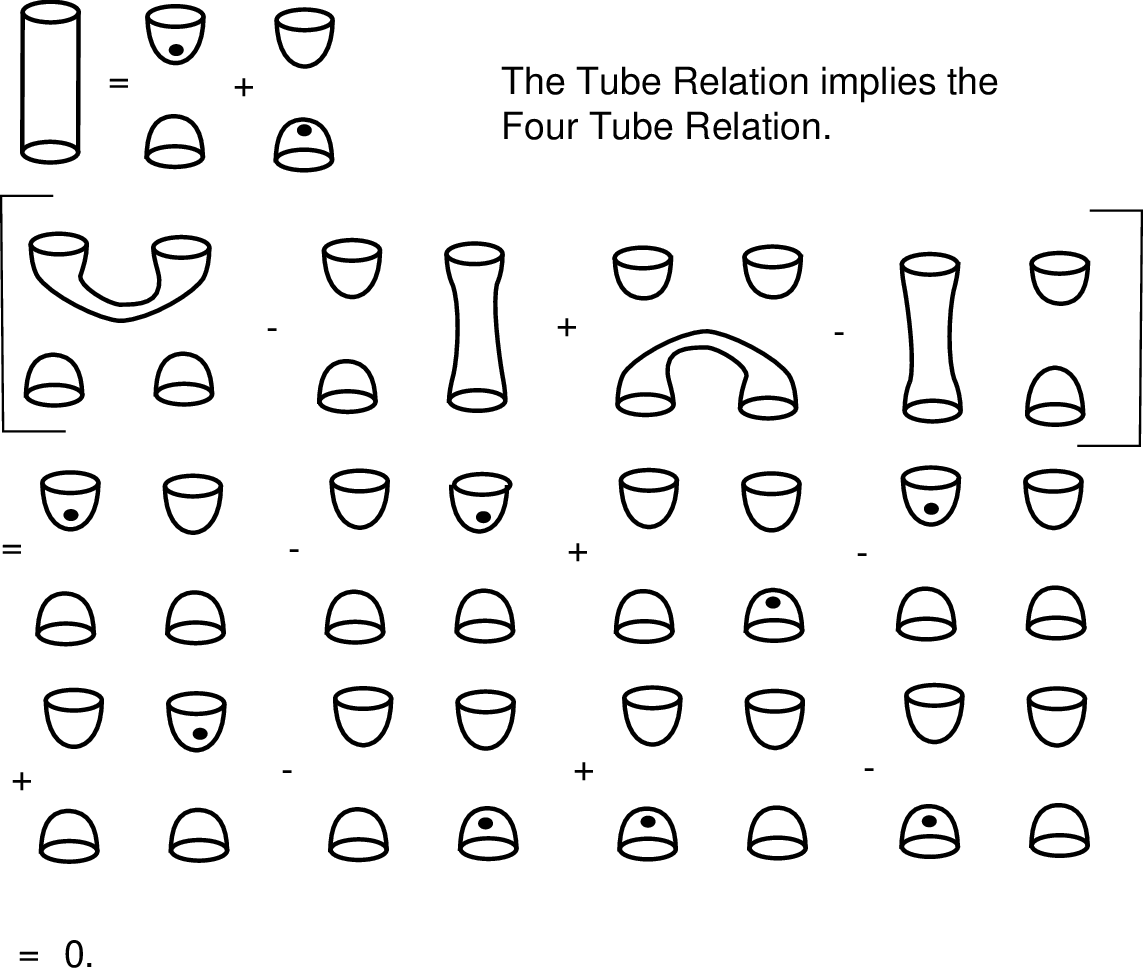}
     \end{tabular}
     \caption{\bf Tube-Cutting Relation Implies Four-Tube Relation}
     \label{tubefourtube}
\end{center}
\end{figure}

\section{\bf Frobenius Algebras Implied by the Tube-Cutting Relation}
In this section we will assume that there is a Frobenius algebra ${\mathcal A}$ that is a ring with identity element $1$ and has an element $x$ that commutes with $1$ and and that $1$ and $x$ are linearly independent over the ring $k.$  We assume that $2$ is an invertible element in the ring $k$. We further assume that the 
dot in the tube-cutting relation stands for the element $x.$ And we assume that the tube-cutting relation is satisfied. As we have seen in Figure 8, this means that 
$$a = \epsilon(ax)1 + \epsilon(a)x$$ for all $a$ in the algebra ${\mathcal A}.$  Thus we shall refer to this equation as the {\it algebraic tube-cutting relation}. At this point we will not make any further
assumptions. As we shall see, these assumptions are sufficient for us to derive a generalization of the Frobenius algebra that works successfully to produce Khovanov Homology. In this way, we see that
the Bar-Natan cobordism picture for the Khovanov invaraiant provides a diagrammatic/topological background from which the basic algebra for the homology can be derived. In another form of exposition, we could have started with only cobordisms and the notion of an abstract complex. Then the particularities of the algebra would be seen as a consequence of the general chain homotopy theory.\\

The approach described above is implicit in Bar--Natan \cite{BN} and it has been carried out in detail by Naot \cite{Naot}. In our work below, we shall restrict ourselves to the consequences of the
tube--cutting relation. Over a general ring, one gets a universal Frobenius algebra $k[x]/(x^{2} - h x)$ with the comultiplication given by $1 \longrightarrow  1 \otimes x + x \otimes 1 - h 1 \otimes 1$ and 
$x  \longrightarrow  x \otimes  x.$ This was also worked out by Naot.\\

Using the algebraic tube-cutting relation, we can write
$$x = \epsilon(x^2)1 + \epsilon(x)x$$
and
$$1 = \epsilon(x)1 + \epsilon(1)x.$$
By linear independence, we conclude that
$$\epsilon(x) = 1, \epsilon(x^2) = 0$$ and $$\epsilon(1) = 0.$$
Furthermore
$$x^2 = \epsilon(x^3)1 + \epsilon(x^2)x,$$ whence
$$x^2 = \epsilon(x^3)1 = t 1$$ where $t  \in k.$
Now look at the coproduct in ${\mathcal A.}$ In Figure 18 we have shown how to expand the cobordism for the coproduct into a sum of terms involving $x, x^2$ and the unit and the counit.
As Figures 19 illustrates, this implies that 
$$\Delta(1) = 1 \otimes x + x \otimes 1$$ and 
$$\Delta(x) = t (1 \otimes 1) + x \otimes x.$$
These equations define a more general Frobenius algebra that can still be used to define a homology theory for knots and links that is invariant under the Reidemeister moves.
Here is a summary of what we have just done.\\

We have produced a Frobenius algebra ${\mathcal A} = k[x]/(x^{2} -t1)$  with $t$ an arbitrary element of the base ring $k,$ and 
$$x^2 = t 1,$$
$$\Delta(1) = 1 \otimes x + x \otimes 1,$$
$$\Delta(x) = t(x \otimes x) + 1 \otimes 1,$$
$$\epsilon(x) = 1,$$
$$\epsilon(1) = 0.$$
For any value of $t$ this algebra satisfies the tube-cutting relation, and so will yield a homology theory that is invariant under the Reidemeister moves. With $t = 0$ we obtain the original 
Frobenius algebra for Khovanov Homology that we have studied in this paper. For $t=1$ we obtain the {\it Lee Homology} that will appear in the next section of this paper.\\ 

\begin{figure}
     \begin{center}
     \begin{tabular}{c}
     \includegraphics[width=9cm]{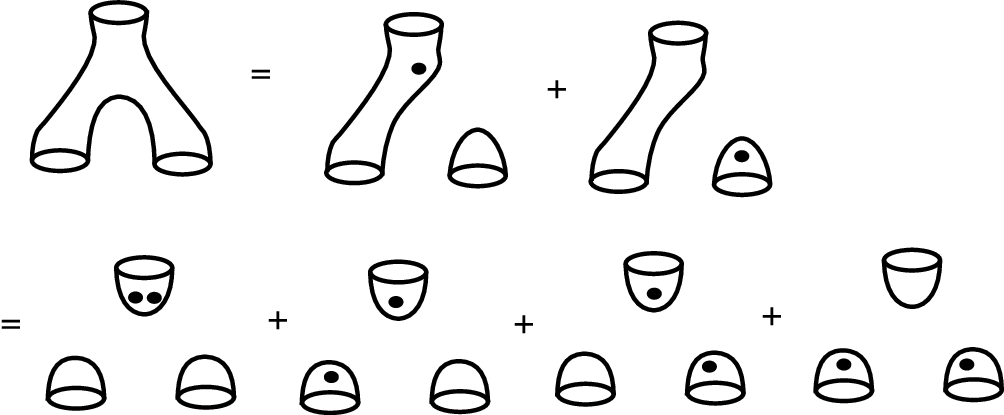}
     \end{tabular}
     \caption{\bf Coproduct Via Tube-Cutting Relation}
     \label{coproduct}
\end{center}
\end{figure}

\begin{figure}
     \begin{center}
     \begin{tabular}{c}
     \includegraphics[width=9cm]{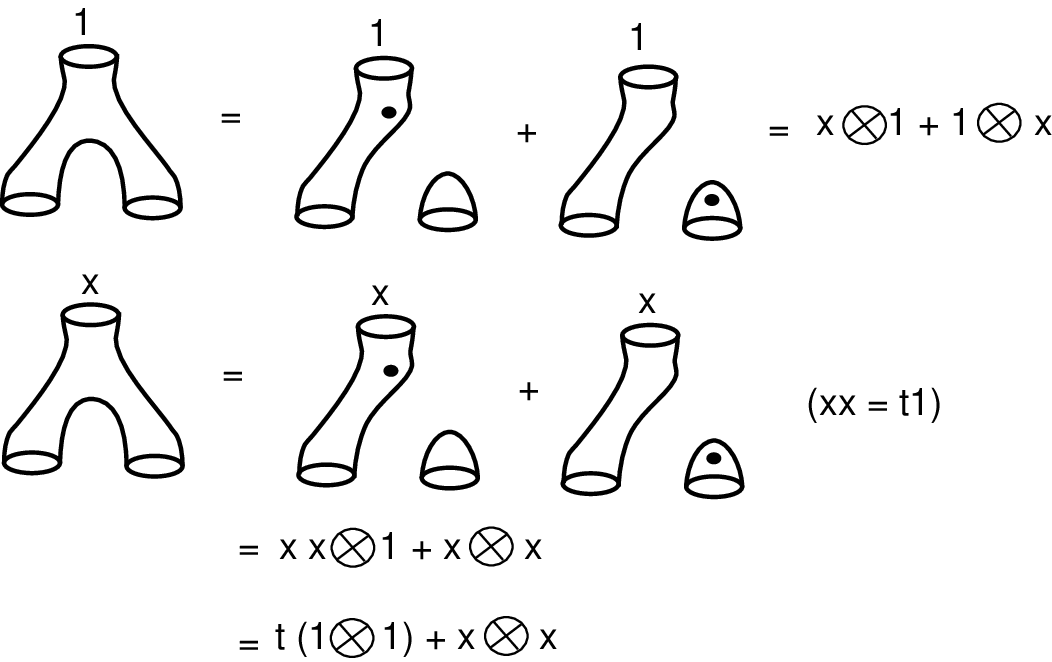}
     \end{tabular}
     \caption{\bf Coproducts of $1$ and $x$ Via Tube-Cutting Relation}
     \label{coproducts}
\end{center}
\end{figure}

\section{\bf Other Frobenius Algebras and Rasmussen's Theorem}
Lee \cite{Lee} makes another homological invariant of knots and links by using a different
Frobenius algebra. She takes the algebra ${\mathcal A} = k[x]/(x^{2} -1)$  with
$$x^2 = 1,$$
$$\Delta(1) = 1 \otimes x + x \otimes 1,$$
$$\Delta(x) = x \otimes x + 1 \otimes 1,$$
$$\epsilon(x) = 1,$$
$$\epsilon(1) = 0.$$
This gives a link homology theory that is distinct from Khovanov homology. In this theory, the quantum
grading $j$ is not preseved, but we do have that $$j(\partial(\alpha)) \ge j(\alpha)$$
for each chain $\alpha$ in the complex. This means that {\em one can use $j$ to filter the chain complex for the Lee homology.} The result is a spectral sequence that starts from Khovanov homology and converges to Lee homology.
\bigbreak

Lee homology is simple. One has that the dimension of the Lee homology is equal to 
$2^{comp(L)}$ where $comp(L)$ denotes the number of components of the link $L.$
Up to homotopy, Lee's homology has a vanishing differential, and the complex behaves well under link
concondance. In his paper \cite{BN1} Dror BarNatan remarks ``In a beautiful article Eun Soo Lee
introduced a second differential on the Khovanov complex of a knot (or link) and showed that 
the resulting (double) complex has non-interesting homology. This is a very interesting result."
Rasmussen \cite{JR} uses Lee's result to define invariants of links that give lower bounds for the 
four-ball genus, and determine it for torus knots. This gives an (elementary) proof of a conjecture
of Milnor that had been previously shown using gauge theory by Kronheimer and Mrowka \cite{KM}.
\bigbreak

Rasmussen's result uses the Lee spectral sequence. We have the quantum ($j$) grading for a 
diagram $K$ and the fact that for Lee's algebra $j(\partial(s)) \ge j(s).$ Rasmussen uses a normalized
version of this grading denoted by $g(s)$.  Then one makes a filtration
$F^{k}C^{*}(K) = \{v \in C^{*}(K) | g(v) \ge k\}$ and given $\alpha \in Lee^{*}(K)$ define
$$S(\alpha) := max\{g(v)| [v] = \alpha \}$$
$$s_{min} (K) := min\{ S(\alpha) | \alpha \in Lee^{*}(K), \alpha \ne 0\}$$
$$s_{max} (K) := max\{ S(\alpha) | \alpha \in Lee^{*}(K), \alpha \ne 0\}$$
and
$$s(K):= (1/2)(s_{min} (K) + s_{max} (K)).$$
This last average of $s_{min}$ and $s_{max}$ is the Rasmussen invariant.
\bigbreak

We now enter the following sequence of facts:
\begin{enumerate}
\item $s(K) \in Z.$
\item $s(K)$ is additive under connected sum.
\item If $K^{*}$ denotes the mirror image of the diagram $K,$ then $$s(K^{*}) = - s(K).$$
\item If $K$ is a positive knot diagram (all positive crossings), then $$s(K) = -r + n +1$$ where 
$r$ denotes the number of loops in the canonical oriented smoothing (this is the same as the 
number of Seifert circuits in the diagram $K$) and $n$ denotes the number of crossings in $K.$
\item For a torus knot $K_{a,b}$ of type $(a,b),$ $s(K_{a,b}) = (a-1)(b-1).$
\item $|s(K)| \le 2 g^{*}(K)$ where $g^{*}(K)$ is the least genus spanning surface for $K$ in the four ball.
\item $g^{*}(K_{a,b}) = (a-1)(b-1)/2.$ This is Milnor's conjecture.
\end{enumerate}

\noindent This completes a very skeletal sketch of the construction and use of Rasmussen's invariant.
\bigbreak

\section  {\bf The Simplicial Structure of Khovanov Homology}
Let ${\mathcal S}$ denote the set of (standard) bracket states for a link diagram $K.$
One way to describe the Khovanov complex is to associate to each state loop $\lambda$
a module $V$ isomorphic to the algebra  $k[x]/(x^2)$ with coproduct as we have described in the previous sections. The generators $1$ and $x$ of this algebra can then be regarded as the two possible
enhancements of the loop $\lambda.$ In the same vein we associate to a state $S$ the tensor
product of copies of $V$, one copy for each loop in the state. The local boundaries are defined exactly as before, and the Khovanov complex is the direct sum of the modules associated with the states of
the link diagram. We will use this point of view in the present section, and we shall describe Khovanov
homology in terms of the  $n$-cube category and an associated simplicial object.
The purpose of this section is to move towards, albeit in an abstract manner, a description of Khovanov
homology as the homology of a topological space whose homotopy type is an invariant of the knot of the underlying knot or link. We do not accomplish this aim, but the constructions
given herein may move toward that goal. An intermediate possibility would be to replace the
Khovanov homology by an abstract space or simplicial object whose generalized homotopy type was an
invariant of the knot or link.  
\bigbreak

Let ${\mathcal D}^{n} = \{A,B \}^{n}$ be the $n$-cube category whose objects are the $n$-sequences from
the set $\{A, B \}$ and whose morphisms are arrows from sequences with greater numbers of $A$'s
to sequences with fewer numbers of $A$'s. Thus ${\mathcal D}^{n}$ is equivalent to the poset category
of subsets of $\{ 1,2,\cdots n \}.$ Let ${\mathcal M }$ be a pointed category with finite sums, and let
${\mathcal F} : {\mathcal D}^{n}  \longrightarrow  {\mathcal M } $ be a functor. In our case ${\mathcal M }$ is a category of
modules (as described above) and ${\mathcal F}$ carries $n$-sequences to certain tensor powers corresponding to the standard bracket states of a knot or link $K.$ We map sequences to states
by choosing to label the crossings of the diagram $K$ from the set  $\{ 1,2,\cdots n \},$ and letting
the functor take abstract $A$'s and $B$'s in the cube category to smoothings at those crossings of type $A$ or type $B.$ Thus each sequence in the cube category is associated with a unique state of $K$
when $K$ has $n$ crossings. Nevertheless, we shall describe the construction more generally.
\bigbreak

For the functor ${\mathcal F}$ we first construct a semisimplicial object ${\mathcal C(F)}$ over ${\mathcal M }$,
where a semisimplicial object is a simplicial object without degeneracies. This means that it
has partial boundaries analogous to the partial boundaries that we have discussed before but none
of the degeneracy maps that are common to simplicial theory (see \cite{May} Chapter 1). For $k \ge 0$
we set $${\mathcal C(F)}_{k} = \oplus_{v  \in {\mathcal D}^{n}_{k}} {\mathcal F}(v)$$
where ${\mathcal D}^{n}_{k}$
denotes those sequences in the cube category with $k$ $A$'s. Note that we are indexing dually
to the upper indexing in the Khovanov homology sections of this paper where we counted the number of
$B$'s in the states.
\bigbreak

We introduce face operators (partial boundaries in our previous terminology)
$$d_{i} : {\mathcal C(F)}_{k} \longrightarrow {\mathcal C(F)}_{k-1}$$ for $0 \le i \le k$ with $k \ge 1$ as
follows: $d_{i}$ is trivial for $i = 0$ and otherwise $d_{i}$ acts on ${\mathcal F}(v)$ by the map
${\mathcal F}(v) \longrightarrow {\mathcal F}(v')$ where $v'$ is the sequence resulting from replacing
the $i$-th $A$ by $B.$ The operators $d_{i}$ satisfy the usual face relations of simplicial theory:
$$d_{i} d_{j} = d_{j-1} d_{i}$$ for $i <  j.$
\bigbreak

We now expand ${\mathcal C}(F)$ to a simplicial object ${\mathcal S}(F)$ over ${\mathcal M }$ by applying freely
degeneracies to the ${\mathcal F}(v)$'s. Thus
$${\mathcal S}(F)_{m} = \oplus_{v \in {\mathcal D}^{n}_{k}, k + t = m} \,\, s_{i_{1}} \cdots  s_{i_{t}} {\mathcal F}(v)$$
where $ m  >  i_{1}  >  \cdots  > i_{t} \ge 0 $ and these degeneracy operators are applied freely modulo the usual (axiomatic) relations among themselves and with the face operators. Then ${\mathcal S}(F)$
has degeneracies via formal application of degeneracy operators to these forms, and has face operators
extending those of ${\mathcal C}({\mathcal F}).$ It is at this point we should remark that in our knot theoretic construction
there is only at this point an opportunity for formal extension of degeneracy operators above the number of crossings in the given knot or link diagram since to make specific degeneracies would involve the
creation of new diagrammatic sites. There is a natural construction of this sort and it can be used to give a simpicial homotopy type for Khovanov homology. See \cite{GK}.
\bigbreak

When the functor ${\mathcal F}: {\mathcal D}^{n} \longrightarrow {\mathcal M}$ goes to an abelian category ${\mathcal M},$
as in our knot theoretic case, we can recover the homology groups via
$$H_{\star} N {\mathcal S}({\mathcal F}) \cong H_{\star} {\mathcal C}({\mathcal F})$$
where $N {\mathcal S}({\mathcal F})$ is the normalized chain complex of ${\mathcal S}({\mathcal F}).$ This completes
the abstract simplicial description of this homology.
\bigbreak

\section{Quantum Comments}
States of a quantum system are represented by unit vectors in a Hilbert space. Quantum processes are unitary transformations applied to these state vectors. In an appropriate basis for the HIlbert space,
each basis vector represesents a possible measurement. If $|\psi \rangle$ is a unit vector, then, upon measurement, one of the basis vectors will appear with probability, the absolute square of its coefficient in $|\psi \rangle.$  One can, in principle, find the trace of a given unitary transformation by instantiating it
in a certain quantum system and making repeated measurements on that system. Such a scheme, in the abstract, is called a quantum algorithm, and in the concrete is called a quantum computer. One  well-known quantum algorithm for determining the trace of a unitary matrix is called the ``Hadamard Test"
\cite{NC}.

In \cite{KQuant1} we consider the Jones polynomial and Khovanov homology in a quantum context.
In this section we give a sketch of these ideas. Recall from Section 2 that we have the following
formula for the Jones polynomial.
$$J_{K} = (-1)^{n_{-}} q^{n_{+} - 2n_{-}} \langle K \rangle.$$
Using the enhanced states formulation of 
Section 2, we form a Hilbert space ${\mathcal H}(K)$ with orthonormal basis the set of enhanced states of $K.$ For the Hilbert space we denote a basis element by $|s \rangle$ where $s$ is an enhanced 
state of the diagram $K.$ Now using $q$ as in Section 2, let $q$ be any point on the unit circle in the
complex plane. Define $U_{K}: {\mathcal H}(K) \longrightarrow {\mathcal H}(K)$ by the formula
$$U_{K} |s \rangle = (-1)^{i(s) + n_{-}}q^{j(s) +n_{+} - 2n_{-}} | s \rangle.$$ Then $U_{K}$ defines a unitary transformation of the Hilbert space and we have that $$J_{K } = Trace(U_{K}).$$ 
\bigbreak

The Hadamard Test
applied to this unitary transformation gives a quantum algorithm for the Jones polynomial. 
This is not the most efficient quantum algorithm for the Jones polynomial. Unitary braid group representions can do better \cite{QCJP,3Strand,Ah1}. But this algorithm has the conceptual advantage of being directly related to Khonavov homology. In particular, let $C^{i,j}$ be the subspace of   ${\mathcal H}(K)$ with basis the set of enhanced states $| s \rangle$ with $i(s) = i$ and $j(s) = j.$ Then ${\mathcal H}(K)$ is
the direct sum of these subspaces and we see that ${\mathcal H}(K)$ is identical with the Khovanov
complex for $K$ with coefficients in the complex numbers. Furthermore, letting 
$\partial: {\mathcal H}(K) \longrightarrow {\mathcal H}(K)$ be the boundary mapping that we have defined for the Khovanov complex, we have $$\partial \circ U_{K} + U_{K} \circ \partial = 0.$$
Thus $U_{K}$ induces a mapping on the Khovanov homology of $K.$ As a linear space, the 
Khovanov homology of $K$, 
$${\mathcal Homology}({\mathcal H}(K)) = Kernel( \partial) /Image (\partial)$$
is also a Hilbert space on which $U_{K}$ acts and for which the trace yields the Jones polynomial.
\bigbreak

If we are given more information about the Khovanov homology as a space, for example if we are given
a basis for $H^{i - n_{-}, j - n_{+} + 2n_{-}}(K)$ for each $i$ and $j$, then we can extend $U$ to act
on $H^{i - n_{-}, j - n_{+} + 2n_{-}}(K)$ as an eigenspace with eigenvalue $ t^{i} q^{j}$ where $q$ and
$t$ are chosen unit complex numbers. Then we have an extended $U'_{K}$ with
$$U'_{K} |\alpha\rangle = t^{i} q^{j} | \alpha \rangle$$ for each $\alpha \in H^{i - n_{-}, j - n_{+} + 2n_{-}}(K).$ With this extension we have that the trace of $U'_{K}$ recovers a specialization of the 
Poincar\'e polynomial (Section 4)  for the Khovanov homology.
$$Trace(U'_{K}) =  \Sigma_{i,j} t^{i} q^{j} dim(H^{i - n_{-}, j - n_{+} + 2n_{-}}(K)) = P_{K}(t, q).$$ 
Thus, in principle, we formulate a quantum algorithm for specializations of the Poincar\`e polynomial
for Khovanov homolgy.
\bigbreak

Placing Khovanov homology in an appropriate quantum mechanical, quantum information theoretic, or quantum field theory context is
a fundamental question that has been considered by a number of people, including Sergei Gukov
\cite{Gukov1,Gukov} and Edward Witten \cite{Witten,Witten1,Witten2}. The constructions discussed here are elementary in nature but we would like to know how they interface with other points of view. In particular, if one thinks of the states in the
state expansion of the bracket polynomial as analogs of the states of a physical system such as the 
Potts model in statistical mechanics, then the loop configuration of a given state corresponds to 
a decomposition of the underlying graph of the statistical mechanics model into regions of constant
spin (where spin designates the local variable in the model). Working with a boundary operator, as we did with the Khovanov chain complex, means taking into account adjacency relations among these 
types of physical states.
\bigbreak

\section{Discussion}
The subject of Khovanov homology is part of the larger subject of categorification in general and other
link homologies in particular.  The term {\it categorification} was coined by Crane and Frenkel in their
paper \cite{CF} speculating on the possibility for invariants of four-manifolds via a categorical generalization of Hopf algebras where all structures are moved up one categorical level. Just such a shift is seen in the Khovanov homology where loops that were once scalars become modules and the original Jones polynomial is seen as a graded Euler characteristic of a homology theory. There is now a complex literature on categorifications of quantum groups (aka Hopf algebras) and relationships of this new form of representation theory with the construction of link homology. For this we refer the reader to the following references \cite{CK1,CK2,KhoLauda1,KhoLauda2,Khovanov,SeidelSmith,Webster,Webster1,Webster2}. It is possible that the vision of Crane and Frenkel for the construction of invariants of four dimensional manifolds will come true.\\

Homotopy and spatial homology theories have been constructed that realize Khovanov homology functorially as homotopy of spectra and homology of spaces. See \cite{ET,ELST,KHT}.\\

Other link homology theories are worth mentioning. In \cite{KhoRoz1,KhoRoz2,KhoRoz3} Khovnaov and Rozansky construct a link
homology theory for specializations of the Homflypt polynomial. Their theory extends integrally to a Khovanov homology theory for virtual knots, but no calculations are known at this writing.
Khovanov homology does extend integrally to virtual knot theory as shown by Manturov in \cite{M1}.
The relationship of the Manturov construction to that of Khovanov and Rozansky is not known at this time.  In \cite{DKM} Dye, Kauffman and Maturov show how to modify mod-2 Khovanov homology to categorify the arrow polynomial for virtual knots. This leads to many new calculations and examples
 \cite{Kaestner,KaestnerKauff}. In \cite{DKK} H. Dye, A. Kaestner and L. H. Kauffman, use a version of Manturov's construction and generalize the Rasmussen invariant to virtual knots and links. They determine the virtual four-ball genus for positive virtual knots.\\
 
 In \cite{MEtAl,MEtAlia} Manolescu, Ozsv\'ath, Szab\'o,  Sarkar and Thurston construct combinatorial link homology based on Floer homology that categorifies the Alexander polynomial. Their techniques are quite different from those explained here for Khovanov homology.  The combinatorial definition should be compared with that of Khovanov homology, but it has a flavor that is different, probably due to the fact that it categorifies a determinant that calculates the Alexander polynomial. This Knot Floer Homology theory is very powerful and can detect the three-dimensional genus of a knot (the least genus of an orientable spanning surface for the knot in three dimensional space). Caprau in \cite{Caprau} has a useful version of the tangle cobordism approach to Khovanov homology and Clark, Morisson and Walker
\cite{CMW} have an oriented tangle cobordism theory that is used to sort out the functoriality of Khovanov homology for knot cobordisms. There is another significant variant of Khovanov homology termed
{\it odd Khovanov homology} \cite{Odd}. Attempts to find other global interpretations of Khovanov homology
have led to very significant lines of research \cite{CK1,CK2,SeidelSmith}, and attempts to find general constructions for 
link homology corresponding to the quantum link invariants coming from quantum groups have led to 
research such as that of Webster \cite{Webster1,Webster2} where we now have theories for such constructions that
use the categorifications of quantum groups for classical Lie algebras.
\bigbreak

There have been three applications of Khovanov homology that are particularly worth mentioning.
One, we have discussed in Section 6, is Rasmussen's use of Khovanov homology \cite{JR} to determine the slice genus of torus knots without using gauge theory. Another is the proof by Kronheimer and Mrowka
\cite{KM} that Khovanov homology detects the unknot. The work of Kronheimer and Mrowka interrelates
Khovanov homology with their theory of knot instanton homology and allows them to apply their
gauge theoretic results to obtain this striking result. A proof that Khovanov homolgy detects the unknot by purely combinatorial topological means is unknown at this writing. By the same token, it is still unknown whether the Jones polynomial detects classical knots. Finally, we mention the work of 
Shumakovitch \cite{Shuma} where, by calculating Khovanov homology, he shows many examples of knots that are topologically slice but are not slice in the differentiable category. Here  Khovanov homology circumvents a previous use of gauge theory but the result still depends on deep results of Freedman showing that classical knots of Alexander polynomial $1$ are topologically slice.
\bigbreak


\begin{thebibliography}{99}
\bibitem{Ah1}
D. Aharonov, V. Jones, Z. Landau, A polynomial quantum algorithm for approximating the Jones polynomial,  Algorithmica 55 (2009), no. 3, 395Ð421, quant-ph/0511096.

\bibitem{BN0}
D. Bar-Natan, On Khovanov's categorification of the Jones polynomial, Algebra Geometry Topology Vol 2. (2002) 337-370, math QA/0201043.

\bibitem{BN}  
D. Bar-Natan (2005), Khovanov's homology for
tangles and cobordisms, {\it Geometry
and Topology}, {\bf 9-33}, pp. 1465-1499. arXiv:mat.GT0410495

\bibitem{BN1} 
D. Bar-Natan and S. Morrison, The Karoubi envelope and Lee's degeneration of Khovanov homology. Algebr. Geom. Topol. 6 (2006), 1459Ð1469. math.GT/0606542.

\bibitem{Caprau}
 C. Caprau, The universal sl(2) cohomology via webs and foams, {\em Topology and its Applications}, Volume 156 (2009), 1684-1702. math.GT. arXiv:0802.2848

\bibitem{CMW}
D. Clark, S. Morrison, K. Walker, Fixing the functoriality of Khovanov homology. {\em Geom. Topol.} 13 (2009), no. 3, 1499Ð1582.

\bibitem{CK1}
S. Cautis, J.  Kamnitzer,  Knot homology via derived categories of coherent sheaves. I. The sl(2)-case. {\em Duke Math. J.} 142 (2008), no. 3, 511Ð588.

\bibitem{CK2}
S. Cautis, J. Kamnitzer, Knot homology via derived categories of coherent sheaves. II. slm case. {\em Invent. Math.} 174 (2008), no. 1, 165Ð232. 


\bibitem{CF}
L. Crane, I. Frenkel, Four-dimensional topological quantum field theory, Hopf categories, and the canonical bases. {\em Topology and physics. J. Math. Phys.} 35 (1994), no. 10, 5136Ð5154.

\bibitem{DKK} H. A. Dye, A. Kaestner, L. H. Kauffman, Khovanov Homology, Lee Homology and a Rasmussen Invariant for Virtual Knots, (submitted to Osaka J. Math.),
arXiv:1409.5088.

\bibitem{ET}
Everitt,Brent;Turner,Paul, The homotopy theory of Khovanov homology, Algebr. Geom. Topol. 14 (2014), no. 5, 2747Ð2781.
arXiv:1112.3460.

\bibitem{KHT}
Lipshitz, Robert; Sarkar, Sucharit A Khovanov stable homotopy type. J. Amer. Math. Soc. 27 (2014), no. 4, 983Ð1042. 

\bibitem{ELST}
 Everitt,Brent; Lipshitz,Robert; Sarkar, Sucharit; Turner,Paul,  Khovanov homotopy types and the Dold-Thom functor. arXiv:1202.1856.
 
 \bibitem{GK} C. Gomes and L. H. Kauffman, An Application of the Dold-Kan Theorem to the Homotopy Theory of Link Homology, (in preparation).
 
\bibitem{Gukov1}
S. Gukov, Gauge theory and knot homologies. Fortschr. Phys. 55 (2007), no. 5-7, 473Ð490.

\bibitem{Gukov} S. Gukov, Surface operators and knot homologies, arXiv:0706.2369.


\bibitem{JO}
V.F.R.  Jones, A polynomial invariant for links via von Neumann algebras,
{\it Bull. Amer. Math. Soc.}, {\bf 129} (1985), 103--112.

\bibitem{JO1}
 V.F.R. Jones.  Hecke algebra representations of braid groups and link polynomials,  {\it Ann. of Math.}  {\bf 126}, (1987), pp. 335-338.

\bibitem{JO2}
V.F.R. Jones, On knot invariants related to some statistical mechanics models, {\it  Pacific J. Math.}, {\bf 137}, no. 2 (1989), pp. 311-334.

\bibitem{KaB}
L.H. Kauffman, State models and the Jones polynomial, {\it Topology}, {\bf 26} (1987),
395--407.

bibitem{KaAM}
L. H. Kauffman. New invariants in knot theory,  {\it Amer. Math. Monthly} {\bf 95},  no. 3, (1988), pp. 195Ð242.

 \bibitem{KA89}
L.H. Kauffman,  Statistical mechanics and the Jones polynomial, {\it  AMS
Contemp. Math. Series},  {\bf 78} (1989), 263--297.

\bibitem {KP}
L.H. Kauffman,  Knots and Physics, World Scientific Publishers (1991),
Second Edition (1993), Third Edition (2002), Fourth Edition (2012).

\bibitem{KQuant}
L H. Kauffman,  Topological quantum information, Khovanov homology and the Jones polynomial. Topology of algebraic varieties and singularities, 245Ð264, {\em Contemp. Math.}, 538, Amer. Math. Soc., Providence, RI, 2011. arXiv:1001.0354.

\bibitem{KTrieste} L. H. Kauffman, Khovanov Homology. in ``Introductory Lectures in 
Knot Theory", K\&E Series Vol. 46, edited by Kauffman, Lambropoulou, Jablan and Przytycki,
World Scientiic 2011, pp. 248 - 280.

\bibitem{KQuant1} L. H. Kauffman. A quantum model for the Jones polynomial, Khovanov
Homology and Generalized Simplicial Homology, {\it AMS
Contemp Math}, {\bf 536}, Cross Disciplinary Advances in Quantum Computing, ed.
by Mahdavi et al., pp. 75 - 94.


\bibitem{QCJP}
L.H. Kauffman, Quantum computing and the
Jones polynomial, in {\it Quantum Computation and Information}, S. Lomonaco,
Jr. (ed.), AMS CONM/305, 2002, pp.~101--137.  math.QA/0105255.

\bibitem{3Strand} L. H.  Kauffman and S.  Lomonaco Jr., A Three-stranded quantum algorithm for the Jones polynonmial, in
{\it Quantum Information and Quantum Computation V}, Proceedings of Spie, April 2007, edited by E.J. Donkor, A.R. Pirich and H.E. Brandt, pp. 65730T1-17, Intl Soc. Opt. Eng.

\bibitem{Kho} M. Khovanov, A categorification of the Jones polynomial, Duke J. Math. 101(200), no. 3, 359-426, mathQA/9908171.

\bibitem{Kock}
J. Kock, ``Frobenius Algebras and 2D Topological Quantum Field Theories", London Mathematical Society, Student Texts No. 59 (2004).

\bibitem{KM}
P. B. Kronheimer, T. S. Mrowka, Khovanov homology is an unknot-detector.  math.GT.arXiv:1005.4346.

\bibitem{Lee} E.  S.  Lee, An endomorphism of the Khovanov invariant, {\em Adv. Math.} 197 (2),
554-586 (2005). math.GT/0210213.

\bibitem{Naot}  G. Naot,  The universal Khovanov link homology theory. Algebr. Geom. Topol. 6 (2006), 1863Ð1892. 

\bibitem{QKnots} S. J.  Lomonaco Jr. and L. H.  Kauffman, Quantum Knots and Mosaics,
{\it Journal of Quantum Information Processing}, Vol. 7, Nos. 2-3, (2008), pp. 85 - 115.  arxiv.org/abs/0805.0339.

\bibitem{DKM} H.A.  Dye, L.H.  Kauffman, V.O.~ Manturov,
On two categorifications of the arrow polynomial for virtual knots, arXiv:0906.3408.

\bibitem{M1} V.O.  Manturov, Khovanov homology for virtual links with arbitrary coefficients,
math.GT/0601152. (Russian) Izv. Ross. Akad. Nauk Ser. Mat. 71 (2007), no. 5, 111--148.
J. Knot Theory Ramifications 16 (2007), no. 3, 345Ð377.

\bibitem{Kaestner} A. Kaestner Ph.D. Thesis, University of Illinois at Chicago (2011).

\bibitem{KaestnerKauff} Kaestner, Aaron M.; Kauffman, Louis H. Parity, skein polynomials and categorification. J. Knot Theory Ramifications 21 (2012), no. 13, 1240011, 56 pp.


\bibitem{KhoLauda1}
M. Khovanov, A. Lauda, A categorification of quantum sl(n). {\em Quantum Topol. }1 (2010), no. 1, 1Ð92.

\bibitem{KhoLauda2}
M. Khovanov, A. Lauda, A diagrammatic approach to categorification of quantum groups II.{\em  Trans. Amer. Math. Soc. } 363 (2011), no. 5, 2685Ð2700.


\bibitem{Khovanov}
M. Khovanov,  Categorifications from planar diagrammatics. {\em Jpn. J. Math.} 5 (2010), no. 2, 153Ð181.

\bibitem{KhoRoz1}
M. Khovanov, L. Rozansky, Virtual crossings, convolutions and a categorification of the SO(2N) Kauffman polynomial. {\em J. Gškova Geom. Topol. GGT }1 (2007), 116Ð214.

\bibitem{KhoRoz2}
M. Khovanov, L. Rozansky, Matrix factorizations and link homology. {\em Fund. Math.} 199 (2008), no. 1, 1Ð91.

\bibitem{KhoRoz3}
M. Khovanov, L. Rozansky, Matrix factorizations and link homology. II. {\em Geom. Topol. }12 (2008), no. 3, 1387Ð1425.

\bibitem{MEtAl}
C. Manolescu, P. Ozsv\'ath, S. Sarkar, A combinatorial description of knot Floer homology. {\em Ann. of Math.} (2) 169 (2009), no. 2, 633Ð660.

\bibitem{MEtAlia}
C. Manolescu, P. Ozsv\'ath, Z. Szab\'o, D. Thurston, On combinatorial link Floer homology. {\em Geom. Topol.} 11 (2007), 2339Ð2412.

\bibitem{Odd}
P. Ozsv\'ath, J. Rasmussen, and Z. Szab\'o, Odd Khovanov homology.{\it  Algebr. Geom. Topol.} {\bf 13} (2013), no. 3, 1465Ð1488.
http://arxiv.org/abs/0710.4300.


\bibitem{Pr} J. H.  Przytycki, When the theories meet:\ Khovanov homology as Hochschild
homology of links, arXiv:math.GT/0509334.

\bibitem{Shuma}
A. N. Shumakovitch,  Rasmussen invariant, slice-Bennequin inequality, and sliceness of knots, {\em J. Knot Theory Ramifications} 16 (2007), no. 10, 1403Ð1412. 

\bibitem{SeidelSmith}
P. Seidel, I. A. Smith,  A link invariant from the symplectic geometry of nilpotent slices. {\em Duke Math. J. }134 (2006), no. 3, 453Ð514.


\bibitem{JR}
J. Rasmussen, Khovanov homology and the slice genus, {\it Invent. Math.} {\bf 182} (2010), no. 2, 419Ð447.
math.GT/0402131.



\bibitem{Lo} J-L.  Loday, Cyclic Homology, Grund. Math. Wissen.
Band 301, Springer-Verlag, Berlin, 1992 (second edition, 1998).


\bibitem{May} J. P.  May, Simplicial Objects in Algebraic Topology, University of Chicago Press
(1967).

\bibitem{NC} M. A.  Nielsen and I. L.  Chuang, Quantum Computation and Quantum Information, Cambridge University Press (2000).


\bibitem {Viro} O.  Viro (2004), Khovanov homology, its definitions and ramifications,
{\it Fund. Math.}, {\bf 184} (2004), pp. 317-342.


\bibitem{Webster}
B. Webster, Khovanov-Rozansky homology via a canopolis formalism. {\em Algebr. Geom. Topol.} 7 (2007), 673Ð699.

\bibitem{Webster1}
 B. Webster, Knot invariants and higher representation theory I: diagrammatic and geometric categorification of tensor products.  math.GT.arXiv:1001.2020.
	
\bibitem{Webster2}
B. Webster, Knot invariants and higher representation theory II: the categorification of quantum knot invariants. math.GT. arXiv:1005.4559

\bibitem{Witten} E. Witten,  Quantum Field Theory and the Jones Polynomial, {\it  Comm. in Math. Phys.},
{\bf 121}, (1989), 351-399. 

	
\bibitem{Witten1}
E. Witten, Knot Invariants from Four-Dimensional Gauge Theory. Davide Gaiotto, Edward Witten. physics.hep-th. arXiv:1106.4789
	
\bibitem{Witten2}
E. Witten, Fivebranes and Knots. Edward Witten. physics.hep-th. arXiv:1101.3216



\bigbreak

\end{thebibliography}
\end{document}